# BASIC PROPERTIES OF NONLINEAR STOCHASTIC SCHRÖDINGER EQUATIONS DRIVEN BY BROWNIAN MOTIONS

By Carlos M. Mora[1] and Rolando Rebolledo[2]

*Universidad de Concepción and Pontificia Universidad Católica de Chile*

The paper is devoted to the study of nonlinear stochastic Schrödinger equations driven by standard cylindrical Brownian motions (NSSEs) arising from the unraveling of quantum master equations. Under the Born–Markov approximations, this class of stochastic evolutions equations on Hilbert spaces provides characterizations of both continuous quantum measurement processes and the evolution of quantum systems. First, we deal with the existence and uniqueness of regular solutions to NSSEs. Second, we provide two general criteria for the existence of regular invariant measures for NSSEs. We apply our results to a forced and damped quantum oscillator.

## 1. Introduction.

1.1. *Objectives.* In this work we focus on the nonlinear stochastic Schrödinger equations driven by standard cylindrical Brownian motions that describe open quantum systems under the Born–Markov approximation. More precisely, this paper concentrates on stochastic evolution equations on a separable complex Hilbert space $(\mathfrak{h}, \langle \cdot, \cdot \rangle)$ of the form

$$(1) \quad X_t = X_0 + \sum_{k=1}^{\infty} \int_0^t (L_k X_s - \operatorname{Re}\langle X_s, L_k X_s \rangle X_s) \, dW_s^k$$
$$+ \int_0^t \left( G X_s + \sum_{k=1}^{\infty} (\operatorname{Re}\langle X_s, L_k X_s \rangle L_k X_s - \tfrac{1}{2} \operatorname{Re}^2 \langle X_s, L_k X_s \rangle X_s) \right) ds,$$

Received October 2006; revised July 2007.
[1]Supported in part by FONDECYT Grants 1040623 and 1070686, PBCT-ACT13 and Grant Milenio ICM P02-049.
[2]Supported in part by FONDECYT Grant 1030552 and PBCT-ACT13.

where $\|X_0\| = 1$ a.s., $W^1, W^2, \ldots$, are real valued independent Wiener processes on a filtered complete probability space $(\Omega, \mathfrak{F}, (\mathfrak{F}_t)_{t \geq 0}, \mathbb{P})$ and $G, L_1, L_2, \ldots$, are linear operators in $\mathfrak{h}$ with $Dom(G) \subset Dom(L_k)$, for any $k \in \mathbb{N}$, such that

$$(2) \qquad 2\mathrm{Re}\langle x, Gx \rangle + \sum_{k=1}^{\infty} \|L_k x\|^2 = 0$$

whenever $x \in Dom(G)$. We are mainly interested in two problems: (i) existence and uniqueness of the solution of (1); and (ii) existence of regular invariant measures for (1).

1.2. *Motivation.* The primary motivation for this article is to develop the mathematical modeling of infinite-dimensional open quantum systems by means of stochastic processes. Our interest in the study of the nonlinear stochastic Schrödinger equation (1) arises mainly from the following three applications.

First, (1) provides characterizations of the evolution of density operators and quantum observables. The classical model of open quantum systems consists of a small quantum system, whose state space is $\mathfrak{h}$ and its internal dynamics is governed by the Hamiltonian $H$, interacting with a heat bath or reservoir. Under the Born–Markov approximation, the time evolution of the density operators (positive operators in $\mathfrak{h}$ with unit trace [16, 56]) obeys the following quantum master equation in Lindblad form:

$$\rho_t = \rho_0 + \int_0^t \left( G\rho_s + \rho_s G^* + \sum_{k=1}^{\infty} L_k \rho_s L_k^* \right) ds.$$

Here $L_1, L_2, \ldots$, take into account the effect of the environment and $G$ is the effective Hamiltonian, that is, $G = -iH - \frac{1}{2}\sum_{k=1}^{\infty} L_k^* L_k$ (see [13, 28] for more details). Let us present a concrete physical model.

EXAMPLE 1 (*Forced and damped quantum oscillator*). Let $\mathfrak{h} = l^2(\mathbb{Z}_+)$. Suppose that $(e_n)_{n \in \mathbb{Z}_+}$ is the canonical orthonormal basis on the space $l^2(\mathbb{Z}_+)$. The linear operators $a^\dagger$ and $a$ are defined in $\{x \in l^2(\mathbb{Z}_+) : \sum_{n \geq 0} n|x_n|^2 < \infty\}$ by $a^\dagger e_n = \sqrt{n+1} e_{n+1}$ and

$$ae_n = \begin{cases} 0, & \text{if } n = 0, \\ \sqrt{n} e_{n-1}, & \text{if } n > 0. \end{cases}$$

Set $N = a^\dagger a$, the number operator.

Choose $H = i\beta_1(a^\dagger - a) + \beta_2 N + \beta_3 (a^\dagger)^2 a^2$. The interaction of the small system with the reservoir is simulated by $L_1 = \alpha_1 a$, $L_2 = \alpha_2 a^\dagger$, $L_3 = \alpha_3 N$, $L_4 = \alpha_4 a^2$, $L_5 = \alpha_5 (a^\dagger)^2$, and $L_6 = \alpha_6 N^2$, where $\beta_1, \beta_2, \beta_3$ are real numbers and $\alpha_1, \ldots, \alpha_6$ are complex numbers. Consider $L_k = 0$ for all $k \geq 7$.

In Example 1, $\mathfrak{h}$ represents the state space of a single mode of a quantized electromagnetic field and the vectors $e_n$, with $n \in \mathbb{Z}_+$, provide the energy levels of the system. Because $a$ destroys a photon, $L_1, L_4$ model photon emissions. The operator $a^\dagger$ describes the creation of a photon (see, e.g., [16, 56]).

Using (1), we can obtain a probabilistic representation of $\rho_t$. Indeed, it is to be expected that $\rho_t = \mathbb{E}|X_t\rangle\langle X_t|$ (see, e.g., [6, 31, 51]). In Dirac notation, $|x\rangle\langle x|$, with $x \in \mathfrak{h}$, stands for the linear operator defined by $|x\rangle\langle x|(z) = \langle x, z\rangle x$ for any $z \in \mathfrak{h}$. Therefore, the probability that a measurement finds the system in the pure state $x$ at time $t \geq 0$ is $\mathbb{E}|\langle x, X_t\rangle|^2$ assuming that $\mathbb{E}|X_0\rangle\langle X_0|$ is the initial density operator and $x$ is a vector of $\mathfrak{h}$ of norm 1. On the other hand, the value of the observable $A$ at time $t$ $\mathcal{T}_t(A)$ satisfies $\mathbb{E}\langle X_t, AX_t\rangle = \mathbb{E}\langle X_0, \mathcal{T}_t(A)X_0\rangle$ (see, e.g., [3, 34, 35]). Recall that quantum observables are represented by self-adjoint operators in $\mathfrak{h}$.

Second, $X_t$ is interpreted as the evolution of the state of a quantum system conditioned on continuous measurement (see, e.g., [4, 8, 19, 30, 51, 59]). For instance, the following example describes the simultaneous monitoring of position and momentum of a quantum system whose evolution is governed by the Hamiltonian $H$ (see, e.g., [32, 57]).

EXAMPLE 2 (*Continuous quantum measurement process*). Let $\mathfrak{h} = L^2(\mathbb{R}, \mathbb{C})$. The position operator $Q: \mathfrak{h} \to \mathfrak{h}$ is given by $Qf(u) = uf(u)$ for any $u \in \mathbb{R}$. The momentum operator $P: \mathfrak{h} \to \mathfrak{h}$ is defined by $P = -iD$, where $Df$ stands for the weak or distributional derivative of a function $f \in H^1(\mathbb{R}, \mathbb{C})$. Then, in the setting of (1) choose $L_1 = \frac{\kappa}{\sigma}Q$ and $L_2 = \kappa\sigma P$, with $\kappa, \sigma \in ]0, +\infty[$. Set $L_k = 0$ for all $k \geq 3$.

Third, (1) plays a relevant role in the numerical simulation of the time evolution of quantum systems (see, e.g., [31, 46, 51, 58]). In fact, using (1), we can overcome the difficulties arising in the direct numerical integration of the master Markov equations (i.e., quantum master equations and Heisenberg equations of motion) in Lindblad form when the dimension of the Hilbert space required for numerical computations is large (see, e.g., [13, 17, 51, 58]).

Another motivation for this paper came from investigations in which (1) represents objective (independent of any observer) trajectories for quantum systems (see, e.g., [31, 51] and the comments of [61]). Furthermore, (1) appears in the quantum filtering theory (see, e.g., [8, 10, 11]). This interesting application has been developed by Belavkin in the framework of quantum stochastic calculus [23, 33, 43, 49].

1.3. *Outline of the paper.* If $\mathfrak{h}$ is finite-dimensional and at most a finite number of $L_k$ are different from 0, then the existence and uniqueness of a strong solution of (1) can be obtained by means of customary techniques employed in stochastic differential equations with locally Lipschitz coefficients (see, e.g., Lemma 5 of [46]). In [5], Barchielli and Holevo covered the existence of a weak solution of (1) when $G$ and $L_1, L_2, \ldots,$ are bounded operators.

Using the linear stochastic Schrödinger equation

$$(3) \qquad \varphi_t(\xi) = \xi + \int_0^t G\varphi_s(\xi)\,ds + \sum_{k=1}^{\infty} \int_0^t L_k \varphi_s(\xi)\,dW_s^k,$$

Girsanov's theorem and Itô's formula, Gatarek and Gisin treated in [29] the existence of solutions of (1) in two examples: (a) $H = 0$, $L_1$ self-adjoint and $L_k = 0$ for any $k \geq 2$; (b) $\mathfrak{h} = L^2(\mathbb{R}, \mathbb{C})$, $H = P^2$, $L_1 = Q$, with $P, Q$ defined as in Example 2, and $L_k = 0$ whenever $k \geq 2$. Exploiting $L_1$ is a self-adjoint operator, [29] also verified the pathwise uniqueness of solutions of (1) in (a) and (b). From the work of Kolokoltsov (see, e.g., [41]) it may deduce a generalization of the results of [29] for some multidimensional versions of the case (b). In [34, 35], Holevo sketched the proof of the existence of a weak solution of (1) in situations where, for any $t > 0$, $(\|\varphi_s(\xi)\|^2)_{s \in [0,t]}$ is uniform integrable and $\mathbb{E}\|\varphi_t(\xi)\|^2 = \mathbb{E}\|\xi\|^2$. To the best of our knowledge, the question of uniqueness of solutions of (1) is still unanswered (even in the case where $G, L_1, L_2, \ldots,$ are bounded operators).

In order to provide a sound basis for the study of infinite-dimensional open quantum systems by means of (1), Section 2 establishes the existence and uniqueness in distribution of the regular solution of (1) under a nonexplosion condition on $G$ and $L_1, L_2, \ldots,$ similar to that used by Chebotarev and Fagnola to prove the Markov property of minimal quantum dynamical semigroups in [14] (see also [15, 23]). To this end, we combine the method of drift transformation (see, e.g., [36, 44, 53]) with the subtleties of (1).

In recent years the large time behavior of quantum Markov semigroups has been the subject of a number of investigations (see, e.g., [2, 24, 25, 26, 60]). However, general results on the large time behavior of stochastic Schrödinger equations on infinite-dimensional Hilbert spaces of the types (1) and (3) are lacking in the literature. It is worth pointing out that Kolokoltsov [40] obtained that $X_t$ is asymptotically similar as $t \to \infty$ to a Gaussian function with certain time-dependent random parameters in case $\mathfrak{h} = L^2(\mathbb{R}, \mathbb{C})$, $H = \hbar P^2/2$, $L_1 = Q$ and $L_2 = L_3 = \cdots = 0$, where $\hbar$ denotes the Planck constant and $P, Q$ are as in Example 2. For this purpose [40] uses an explicit solution of (3) (see, e.g., [9]).

Since $\|X_t\| = 1$ a.s., standard techniques of finite-dimensional stochastic processes yield the existence of stationary distributions for (1) when the dimension of $\mathfrak{h}$ is finite (see, e.g., [22]). A few attempts have been made to develop sufficient conditions for the uniqueness of an invariant measure for (1) in case $\dim \mathfrak{h} < \infty$ (see, e.g., [7, 39]).

Section 3 deals with the basic problem (ii). Under underlying assumptions similar to those in [24], Section 3 states the existence of an invariant probability measure $\Gamma$ for (1) satisfying $\int_{\mathfrak{h}} \|Az\|^2 \Gamma(dz) < \infty$, with $A$ belonging to certain class of linear operators. This regularity property leads to the existence of a quantum stationary state $\rho_\infty$ for which the trace of $A\rho_\infty$ is well defined for a broad class of unbounded operators $A$. Here, our main criterion for the existence of regular invariant measures is based on characteristics of the operators $G$ and $L_k$. Moreover, this criterion involves the existence of a Lyapunov function inherent in the open quantum systems set-up. As a by-product, Section 3 provides the continuous dependence of the distribution of $X$ on the initial data and the Markov property of $X$.

Finally, Section 4 illustrates our main results with a forced and damped quantum oscillator. We select Example 1 as a model problem due to the role played by the one-dimensional quantum harmonic oscillators in the understanding of quantum systems (see, e.g., [16, 56]).

**2. Existence and uniqueness.** This section provides a detailed study of the problem (i) in a context similar to that of [14] (see also [15, 23, 35]). We begin by specifying notation.

2.1. *Notation.* Throughout this paper, the scalar product $\langle \cdot, \cdot \rangle$ is linear in the second variable and anti-linear in the first one. Furthermore, we assume that $(e_n)_{n \in \mathbb{Z}_+}$ is an orthonormal basis of $\mathfrak{h}$. Let $\mathfrak{h}_n$ be the linear manifold spanned by $e_0, \ldots, e_n$. We define $P_n : \mathfrak{h} \to \mathfrak{h}_n$ to be the orthogonal projection of $\mathfrak{h}$ over $\mathfrak{h}_n$.

Let $A$ be a linear operator in $\mathfrak{h}$. Then $Dom(A)$ stands for the domain of $A$ and $A^*$ denotes the adjoint of $A$. The function $\pi_A : \mathfrak{h} \to \mathfrak{h}$ is defined by

$$\pi_A(x) = \begin{cases} x, & \text{if } x \in Dom(A), \\ 0, & \text{if } x \notin Dom(A). \end{cases}$$

Suppose that $C$ is a self-adjoint positive operator in $\mathfrak{h}$. For any $x, y \in Dom(C)$, we set $\langle x, y \rangle_C = \langle x, y \rangle + \langle Cx, Cy \rangle$ and $\|x\|_C = \sqrt{\langle x, x \rangle_C}$. Since $C$ is a closed operator, $(Dom(C), \langle \cdot, \cdot \rangle_C)$ is a Hilbert space. We write $A$ instead of $A \circ \pi_C$ as soon as the context avoid any confusion.

The Borel $\sigma$-algebra of the topological space $E$ is written $\mathfrak{B}(E)$.

2.2. *Main result.* In this section we suppose the existence of a reference operator $C$ with the following properties.

HYPOTHESIS 1. The linear operator $C : \mathfrak{h} \to \mathfrak{h}$ is a self-adjoint positive operator such that:

(H1.1) $Dom(C) \subset Dom(G) \cap Dom(G^*)$.

(H1.2) For any $n \in \mathbb{Z}_+$, $\mathfrak{h}_n \subset Dom(C)$ and $\sum_{k=1}^{\infty} \|L_k^* e_n\|^2 < \infty$.

(H1.3) There exist constants $\alpha, \beta \in [0, +\infty[$ satisfying

$$2\mathrm{Re}\langle Cx, CP_n Gx \rangle + \sum_{k=1}^{\infty} \|CP_n L_k x\|^2 \leq \alpha(\|Cx\|^2 + \|x\|^2 + \beta)$$

for any $n \in \mathbb{Z}_+$ and $x \in \mathfrak{h}_n$.

(H1.4) For all $x \in Dom(C)$, $\sup_{n \in \mathbb{Z}_+} \|CP_n x\| \leq \|Cx\|$.

REMARK 1. Hypothesis 1 is a nonexplosion condition inherent in the open quantum systems context (see, e.g., [15, 23]). It applies to a broad range of applications as, for instance, models for heavy ion collisions [14], quantum oscillators (see, e.g., Section 4.3 of [23]) and quantum exclusion processes [47].

The notion of a regular solution of (1) given below is deeply inspired on the concept of a smooth classical solution of a partial differential equation. Loosely speaking, Definition 1 replaces the partial derivatives of a function by operators $C$ satisfying Hypothesis 1 in order to describe the smoothness of a solution of (3).

DEFINITION 1. Let $C$ satisfy Hypothesis 1. Suppose that $\mathbb{T}$ is either $[0, +\infty[$ or $[0, T]$, provided $T \in [0, +\infty[$. We say that $(\Omega, \mathfrak{F}, (\mathfrak{F}_t)_{t \in \mathbb{T}}, \mathbb{P},$ $(X_t)_{t \in \mathbb{T}}, (W_t^k)_{t \in \mathbb{T}}^{k \in \mathbb{N}})$ is a solution of class $C$ of (1) with initial distribution $\theta$ on the interval $\mathbb{T}$ if:

- $W^1, W^2, \ldots$, are real valued independent Brownian motions on the filtered complete probability space $(\Omega, \mathfrak{F}, (\mathfrak{F}_t)_{t \in \mathbb{T}}, \mathbb{P})$.
- $(X_t)_{t \in \mathbb{T}}$ is an $\mathfrak{h}$-valued process with continuous sample paths such that the law of $X_0$ coincides with $\theta$ and $\mathbb{P}(\|X_t\| = 1$ for all $t \in \mathbb{T}) = 1$.
- For every $t \in \mathbb{T}$, $X_t \in Dom(C)$ $\mathbb{P}$-a.s. and $\sup_{s \in [0,t]} \mathbb{E}_\mathbb{P} \|CX_s\|^2 < \infty$.
- $\mathbb{P}$-a.s., for all $t \in \mathbb{T}$,

(4) $$X_t = X_0 + \int_0^t G(X_s)\,ds + \sum_{k=1}^{\infty} \int_0^t L_k(X_s)\,dW_s^k,$$

where, for any $y \in \mathfrak{h}$,

(5) $$\begin{aligned} G(y) &= G \circ \pi_C(y) \\ &\quad + \sum_{k=1}^{\infty} (\mathrm{Re}\langle y, L_k \circ \pi_C(y)\rangle L_k \circ \pi_C(y) - \tfrac{1}{2}\mathrm{Re}^2 \langle y, L_k \circ \pi_C(y)\rangle y) \end{aligned}$$

and for any $y \in \mathfrak{h}$ and $k \in \mathbb{N}$,

(6) $$L_k(y) = L_k \circ \pi_C(y) - \mathrm{Re}\langle y, L_k \circ \pi_C(y)\rangle y.$$

For abbreviation, we simply say $(\mathbb{P}, (X_t)_{t \in \mathbb{T}}, (W_t)_{t \in \mathbb{T}})$ is a $C$-solution of (1) when no confusion can arise.

The following theorem asserts the existence and uniqueness of the weak (in the probabilistic sense) regular solution of (1).

THEOREM 1. *Let $C$ satisfy Hypothesis 1. Suppose that $\theta$ is a probability measure on $\mathfrak{B}(\mathfrak{h})$ such that $\theta(Dom(C) \cap \{x \in \mathfrak{h} : \|x\| = 1\}) = 1$ and $\int_{\mathfrak{h}} \|Cx\|^2 \theta(dx) < \infty$. Then, (1) has a unique $C$-solution $(\mathbb{Q}, (X_t)_{t \geq 0}, (B_t)_{t \geq 0})$ with initial distribution $\theta$. Furthermore, for any $t \geq 0$, we have*

(7) $\quad \mathbb{E}_{\mathbb{Q}} \|CX_t\|^2 \leq \exp(\alpha t)(\mathbb{E}_{\mathbb{Q}} \|CX_0\|^2 + t\alpha(\mathbb{E}_{\mathbb{Q}} \|X_0\|^2 + \beta)).$

Theorem 1 shows that Example 2 equipped with a standard Hamiltonian is mathematically sound.

COROLLARY 1. *Let assumptions of Example 2 hold. Select $H = \alpha P^2 + \beta Q^2$, with $\alpha \geq 0$ and $\beta \in \mathbb{R}$, and consider the self-adjoint operator [in $L^2(\mathbb{R}, \mathbb{C})$] $N = (Q^2 + P^2 - I)/2$. Suppose that $\theta$ is a probability measure on $\mathfrak{B}(L^2(\mathbb{R}, \mathbb{C}))$ such that, for a given $p \in \mathbb{N}$, $\theta(Dom(N^p) \cap \{x \in L^2(\mathbb{R}, \mathbb{C}) : \|x\| = 1\}) = 1$ and $\int_{\mathfrak{h}} \|N^p x\|^2 \theta(dx) < \infty$. Then, (1) has a unique solution of class $N^p$ with initial distribution $\theta$.*

PROOF. According to Subsection 4.1 of [47], we have $N^p$ satisfies Hypothesis 1. Then, Theorem 1 yields the desired result. □

2.3. *Proof of Theorem 1.* Motivated by the nonlinear filtering theory, the proof of Theorem 1 combines characteristics of (1) with classical techniques based on Girsanov's theorem. We start by recalling the result on the existence and uniqueness of the regular strong solution of (3) given by [47] (see also [45]).

DEFINITION 2. *Let $C$ satisfy Hypothesis 1. Suppose that $\mathbb{T}$ is either $[0, +\infty[$ or $[0, T]$ with $T \in [0, +\infty[$. We say that the stochastic process $(\varphi_t(\xi))_{t \in \mathbb{T}}$ is a strong solution of class $C$ of (3) on the interval $\mathbb{T}$ (for simplicity, $C$-strong solution) if:*

- $(\varphi_t(\xi))_{t \in \mathbb{T}}$ is an adapted process taking values in $\mathfrak{h}$ with continuous sample paths.
- For any $t \in \mathbb{T}$, $\mathbb{E}\|\varphi_t(\xi)\|^2 \leq \mathbb{E}\|\xi\|^2$, $\varphi_t(\xi) \in Dom(C)$ $\mathbb{P}$-a.s. and $\sup_{s \in [0,t]} \mathbb{E}\|C\varphi_s(\xi)\|^2 < \infty$.
- $\mathbb{P}$-a.s., for all $t \in \mathbb{T}$,

(8) $\quad \varphi_t(\xi) = \xi + \int_0^t G \circ \pi_C(\varphi_s(\xi)) \, ds + \sum_{k=1}^{\infty} \int_0^t L_k \circ \pi_C(\varphi_s(\xi)) \, dW_s^k.$

THEOREM 2. *Let $C$ satisfy Hypothesis 1. Suppose that $\xi$ is a $\mathfrak{F}_0$-random variable taking values in $\mathfrak{h}$ such that $\xi \in Dom(C)$ a.s. and $\mathbb{E}\|\xi\|_C^2 < \infty$. Assume that $\mathbb{T}$ is either $[0, +\infty[$ or $[0, T]$ whenever $T \in [0, +\infty[$. Then, there exists a unique $C$-strong solution $(\varphi_t(\xi))_{t \in \mathbb{T}}$ of (3). In addition, for all $t \in \mathbb{T}$, we have the following:*

(i) $\mathbb{E}\|\varphi_t(\xi)\|^2 = \mathbb{E}\|\xi\|^2$.
(ii) $\mathbb{E}\|C\varphi_t(\xi)\|^2 \leq \exp(\alpha t)(\mathbb{E}\|C\xi\|^2 + \alpha t(\mathbb{E}\|\xi\|^2 + \beta))$.

We now point out some immediate consequences of (H1.1), (H1.2) and (H1.4) (see [47] for more details).

REMARK 2. Let assumptions (H1.1) and (H1.2) hold. Using the closed graph theorem, we see that $G$ can be considered as a bounded operator from $(Dom(C), \langle \cdot, \cdot \rangle_C)$ into $\mathfrak{h}$. By (2), for any $x \in Dom(C)$, $\sum_{k=1}^{\infty} \|L_k x\|^2 \leq K\|Cx\|_C^2$, with $K > 0$.

REMARK 3. Suppose that conditions (H1.2) and (H1.4) hold. Then, for any $x$ in $Dom(C)$, $\lim_{n \to \infty} CP_n x = Cx$. It follows that

$$Dom(C) = \{x \in \mathfrak{h} : (CP_n x)_{n \in \mathbb{Z}_+} \text{ is a convergent sequence}\},$$

because $C$ is a closed operator.

REMARK 4. Let $C$ satisfy Hypothesis 1. Applying Remarks 2 and 3, we see that $G \circ \pi_C$ and $L_k \circ \pi_C$, with $k \in \mathbb{N}$, are $\mathfrak{B}(\mathfrak{h})/\mathfrak{B}(\mathfrak{h})$-measurable functions.

The following proposition constructs a $C$-solution of (1) on $[0, T]$ whenever $T \in [0, +\infty[$ with the help of (3).

PROPOSITION 1. *Suppose that hypotheses of Theorem 1 hold. Let $(\varphi_t(\xi))_{t \geq 0}$ be the $C$-strong solution of (3), where $\xi$ is distributed according to $\theta$. Define $\mathbb{Q} = \|\varphi_T(\xi)\|^2 \cdot \mathbb{P}$, where $T \in ]0, +\infty[$. For any $t \in [0, T]$, we set*

$$X_t = \begin{cases} \varphi_t(\xi)/\|\varphi_t(\xi)\|, & \text{if } \varphi_t(\xi) \neq 0, \\ 0, & \text{if } \varphi_t(\xi) = 0, \end{cases}$$

*and*

(9) $$B_t^k = W_t^k - \int_0^t \frac{1}{\|\varphi_s(\xi)\|^2} d[W^k, \varphi(\xi)]_s,$$

*with $k \in \mathbb{N}$. Then $(\Omega, \mathfrak{F}, (\mathfrak{F}_t)_{t \in [0,T]}, \mathbb{Q}, (X_t)_{t \in [0,T]}, (B_t^k)_{t \in [0,T]}^{k \in \mathbb{N}})$ is a $C$-solution of (1) with initial distribution $\theta$.*

PROOF. Applying Itô's formula (or Lemma 2.1 of [47]), we obtain

$$\|\varphi_t(\xi)\|^2 = \|\xi\|^2 + \sum_{k=1}^{\infty} \int_0^t 2\mathrm{Re}\langle \varphi_s(\xi), L_k\varphi_s(\xi)\rangle \, dW_s^k. \tag{10}$$

Consider the stopping time $T_n = \inf\{t \geq 0 : \|\varphi_t(\xi)\| > n\} \wedge T$, with $n \in \mathbb{N}$. From Remark 2 and Theorem 2 we see that

$$\sum_{k=1}^{\infty} \mathbb{E} \int_0^{T_n} \mathrm{Re}^2 \langle \varphi_s(\xi), L_k\varphi_s(\xi)\rangle \, ds \leq n^2 K_T(\mathbb{E}\|\xi\|_C^2 + 1),$$

where $K_T$ is a constant depending on $T$. It follows that

$$\left( \sum_{k=1}^{\infty} \int_0^{t \wedge T_n} 2\mathrm{Re}\langle \varphi_s(\xi), L_k\varphi_s(\xi)\rangle \, dW_s^k \right)_{t \in [0,T]}$$

is a square integrable martingale. Conditional Fatou's lemma now shows that $(\|\varphi_t(\xi)\|^2)_{t \in [0,T]}$ is a supermartingale. Since $\mathbb{E}\|\varphi_t(\xi)\|^2 = \mathbb{E}\|\xi\|^2$ for all $t \geq 0$, $(\|\varphi_t(\xi)\|^2)_{t \in [0,T]}$ is a martingale on $(\Omega, \mathfrak{F}, (\mathfrak{F}_t)_{t \in [0,T]}, \mathbb{P})$.

Let $\mathbb{Q}$ be the probability measure on $\mathfrak{F}$ given by $\mathbb{Q} = \|\varphi_T(\xi)\|^2 \cdot \mathbb{P}$, that is, $\mathbb{Q}$ is absolutely continuous with respect to $\mathbb{P}$ and the Radon–Nikodym derivative of $\mathbb{Q}$ with respect to $\mathbb{P}$ is $\|\varphi_T(\xi)\|^2$. The Girsanov–Meyer theorem shows that $B^k$ given by (9) is a continuous $(\Omega, \mathfrak{F}, (\mathfrak{F}_t)_{t \in [0,T]}, \mathbb{Q})$-local martingale. Since

$$B_t^k = W_t^k - \int_0^t \frac{2\mathrm{Re}\langle \varphi_s(\xi), L_k\varphi_s(\xi)\rangle}{\|\varphi_s(\xi)\|^2} \, ds$$

$$\text{for all } t \in [0,T], [B^k, B^j]_t = \begin{cases} 1, & \text{if } k = j, \\ 0, & \text{if } k \neq j. \end{cases} \tag{11}$$

According to Lévy's theorem, $B^1, B^2, \ldots$, are independent Brownian motions on $(\Omega, \mathfrak{F}, (\mathfrak{F}_t)_{t \in [0,T]}, \mathbb{Q})$.

For all $t \in [0,T]$, $\mathbb{Q}(\varphi_t(\xi) = 0) = 0$. Hence, combining (10) with (11) yields

$$\|\varphi_t(\xi)\|^2 = \|\xi\|^2 + \sum_{k=1}^{\infty} \int_0^t \|\varphi_s(\xi)\|^2 4\mathrm{Re}^2\langle X_s, L_k X_s\rangle \, ds$$

$$+ \sum_{k=1}^{\infty} \int_0^t \|\varphi_s(\xi)\|^2 2\mathrm{Re}\langle X_s, L_k X_s\rangle \, dB_s^k. \tag{12}$$

Let $M_t = \sum_{k=1}^{\infty} \int_0^t 2\mathrm{Re}\langle X_s, L_k X_s\rangle \, dB_s^k$. From Remark 2, Theorem 2 and

$$\mathbb{E}_{\mathbb{Q}} \mathrm{Re}^2\langle X_s, L_k X_s\rangle \leq \mathbb{E}_{\mathbb{Q}} \|L_k X_s\|^2$$
$$= \mathbb{E}_{\mathbb{P}} \|L_k \varphi_s(\xi)\|^2,$$

we have $(M_t)_{t \in [0,T]}$ is a continuous square integrable martingale. By (12),

$$\|\varphi_t(\xi)\|^2 = \|\xi\|^2 + \int_0^t \|\varphi_s(\xi)\|^2 \, d(M + [M,M])_s.$$

Therefore, $\|\varphi_t(\xi)\|^2 = \exp(M_t + [M,M]_t/2)\|\xi\|^2$, which implies

$$\|\varphi_t(\xi)\|^{-1} = \exp(-M_t/2 - [M,M]_t/4)\|\xi\|^{-1}.$$

Hence,

$$\begin{aligned}
\|\varphi_t(\xi)\|^{-1} = \|\xi\|^{-1} &- \tfrac{1}{2} \sum_{k=1}^{\infty} \int_0^t \|\varphi_s(\xi)\|^{-1} \operatorname{Re}^2 \langle X_s, L_k \circ X_s \rangle \, ds \\
&- \sum_{k=1}^{\infty} \int_0^t \|\varphi_s(\xi)\|^{-1} \operatorname{Re} \langle X_s, L_k X_s \rangle \, dB_s^k.
\end{aligned} \quad (13)$$

Using Itô's formula, (8) and (13), we obtain (4) with $W$ replaced by $B$. □

REMARK 5. Let assumptions of Proposition 1 hold. By the previous proof, the process $(\|\varphi_t(\xi)\|^2)_{t \geq 0}$ is a martingale on $(\Omega, \mathfrak{F}, (\mathfrak{F}_t)_{t \in [0,T]}, \mathbb{P})$. Then for every $T > 0$, $(\|\varphi_t(\xi)\|^2)_{t \in [0,T]}$ is a uniformly integrable martingale as both Section 4 of [34] and Section 4 of [35] stated without proof.

NOTATION 1. Let $E$ be a normed space. In the sequel, $C([0,T], E)$ stands for the space of all continuous functions from the interval $[0,T]$ to $E$ endowed with the uniform norm $\|\cdot\|_{\infty}$. Moreover, we define $C([0,T], \mathbb{R}^{\infty})$ to be the Cartesian product space $\prod_{k=1}^{\infty} C([0,T], \mathbb{R})$ equipped with the metric

$$d((f^k)_{k \in \mathbb{N}}, (g^k)_{k \in \mathbb{N}}) = \sum_{k=1}^{\infty} 2^{-k} \min\{1, \|f^k - g^k\|_{\infty}\}.$$

Proposition 2 deals with the uniqueness of solutions of class $C$ for (1) in the sense of joint probability law when $\mathbb{T}$ is a bounded interval.

PROPOSITION 2. *Let hypotheses of Theorem 1 hold. Take $(X_t, (B_t^k)^{k \in \mathbb{N}})_{t \in [0,T]}$ as in Proposition 1. If $(\widetilde{\mathbb{Q}}, (\tilde{X}_t)_{t \in [0,T]}, (\tilde{B}_t)_{t \in [0,T]})$ is a C-solution of (1) with initial distribution $\theta$, then $(\tilde{X}_t, (\tilde{B}_t^k)^{k \in \mathbb{N}})_{t \in [0,T]}$ and $(X_t, (B_t^k)^{k \in \mathbb{N}})_{t \in [0,T]}$ have the same finite-dimensional distributions.*

PROOF. We begin by recalling that

$$\mathfrak{B}(C([0,T], \mathbb{R}^{\infty})) = \bigotimes_{k=1}^{\infty} \mathfrak{B}(C([0,T], \mathbb{R})).$$

Then, $(B_t^k)_{t\in[0,T]}^{k\in\mathbb{N}}$ and $(\tilde{B}_t^k)_{t\in[0,T]}^{k\in\mathbb{N}}$ are random variables taking values in $C([0,T],\mathbb{R}^\infty)$. Since $C([0,T],\mathbb{R}^\infty)$ is a Polish space, classical arguments (see, e.g., the proof of Proposition IX.1.4 of [53]) reduce our proof to the case $X_0 = x$ and $\tilde{X}_0 = x$ for $x \in Dom(C)$ with $\|x\| = 1$.

To deal with the above situation, we consider $\Theta = C([0,T],\mathfrak{h}) \times C([0,T],\mathfrak{h}) \times C([0,T],\mathbb{R}^\infty)$ equipped with the usual product topology. For any $f = (f^1, f^2, f^3) \in \Theta$, we set $Z^k(f) = f^k$, provided $k = 1, 2, 3$. By the Yamada–Watanabe construction (see, e.g., proof of Theorem IV.1.1 of [36] or proof of Theorem IX.1.7 of [53]), there exists a probability measure $\mu$ on $\mathfrak{B}(\Theta)$ with the following properties:

- $\mu \circ (Z^1, Z^3)^{-1} = \widetilde{\mathbb{Q}} \circ (X, B)^{-1}$ and $\mu \circ (Z^2, Z^3)^{-1} = \widetilde{\mathbb{Q}} \circ (\tilde{X}, \tilde{B})^{-1}$.
- $(Z_t^3)_{t\in[0,T]}$ is a sequence of real valued independent Brownian motions on $(\Theta, \mathfrak{G}, (\mathfrak{G}_t)_{t\in[0,T]}, \mu)$. Here $(\Theta, \mathfrak{G}, \mu)$ is the completion of $(\Theta, \mathfrak{B}(\Theta), \mu)$ and for each $t \in [0,T]$,

$$\mathfrak{G}_t = \bigcap_{\varepsilon>0} (\sigma(f(s) : s \in [0, (t+\varepsilon) \wedge T]) \cup \{A \in \mathfrak{G} : \mu(A) = 0\}).$$

In the sequel, $j$ is either 1 or 2. From Remark 2 we have

$$\sum_{k=1}^\infty \int_0^T \mathbb{E}_\mu \text{Re}^2 \langle Z_s^j, L_k Z_s^j \rangle \, ds \leq K \int_0^T (1 + \mathbb{E}_\mu \|C Z_s^j\|^2) \, ds,$$

and so $\sum_{k=1}^\infty \int_0^T \mathbb{E}_\mu \text{Re}^2 \langle Z_s^j, L_k Z_s^j \rangle \, ds < \infty$. Then $(H_t)_{t\in[0,T]}$ described by

$$H_t = -\sum_{k=1}^\infty \int_0^t 2\text{Re}\langle Z_s^j, L_k Z_s^j \rangle \, d(Z^3)_s^k$$

is a continuous square integrable martingale. For each $n \in \mathbb{N}$, we define

$$T_n^j = \inf\left\{ t \in [0,T] : \int_0^t \|C Z_s^j\|^2 \, ds > n \right\} \wedge T.$$

Therefore, $T_n^j$ is a stopping time and $[H^{T_n^j}, H^{T_n^j}]_t \leq Km(n+T)$. According to Novikov's criterion (see, e.g., Proposition VIII.1.15 of [53]), we have $(\exp(H_t^{T_n^j} - [H^{T_n^j}, H^{T_n^j}]_t/2))_{t\in[0,T]}$ is a uniformly integrable martingale. For any $k \in \mathbb{N}$, we choose

(14) $$W_t^{k,j} = (Z^3)_t^k + \int_0^{t \wedge T_n^j} 2\text{Re}\langle Z_s^j, L_k Z_s^j \rangle \, ds.$$

Applying Girsanov's theorem, we deduce that $(W_t^{k,j})_{t\in[0,T]}^{k\in\mathbb{N}}$ is a sequence of real valued independent Brownian motions on $(\Theta, \mathfrak{G}, (\mathfrak{G}_t)_{t\in[0,T]}, \mu^j)$, where

$$\mu^j = \exp(H_T^{T_n^j} - [H^{T_n^j}, H^{T_n^j}]_T/2) \cdot \mu.$$

Because $(Z_t^j, Z_t^3)_{t\in[0,T]}$ satisfies (4) with $(X, B)$ replaced by $(Z^j, Z^3)$ (see, e.g., Theorems 8.3 and 8.6 of [48] or Exercise IV.5.16 of [53]), we see that

$$(Z^j)_t^{T_n^j} = x + \int_0^{t\wedge T_n^j} \left( GZ_s^j + \tfrac{3}{2} \sum_{k=1}^\infty \mathrm{Re}^2 \langle Z_s^j, L_k Z_s^j \rangle Z_s^j \right) ds$$

$$- \int_0^{t\wedge T_n^j} \left( \sum_{k=1}^\infty \mathrm{Re}\langle Z_s^j, L_k Z_s^j \rangle L_k Z_s^j \right) ds$$

$$+ \sum_{k=1}^\infty \int_0^{t\wedge T_n^j} L_k(Z_s^j)\, dW_s^{k,j}.$$

For any $t \in [0, T]$, we set $\eta_t = \exp(-H_t/2 + [H, H]_t/4)$. By, for instance, Theorem II.37 of [50],

$$\eta_t^{T_n^j} = 1 - \tfrac{1}{2} \sum_{k=1}^\infty \int_0^{t\wedge T_n^j} \eta_s \mathrm{Re}^2 \langle Z_s^j, L_k Z_s^j \rangle\, ds$$

$$+ \sum_{k=1}^\infty \int_0^{t\wedge T_n^j} \eta_s \mathrm{Re}\langle Z_s^j, L_k Z_s^j \rangle\, dW_s^{k,j}.$$

The Itô formula leads to

$$(Z^j \eta)_t^{T_n^j} = x + \int_0^{t\wedge T_n^j} G(Z^j \eta)_s^{T_n^j}\, ds + \sum_{k=1}^\infty \int_0^{t\wedge T_n^j} L_k(Z^j \eta)_s^{T_n^j}\, dW_s^{k,j}.$$

Let $(\varphi_t^j(x))_{t\in[0,T]}$ be the solution of (8) described by Theorem 2 when the underlying filtered probability space is $(\Theta, \mathfrak{G}, (\mathfrak{G}_t)_{t\in[0,T]}, \mu^j)$. In order to prove that $(Z^j \eta_{t\wedge T_n^j})_{t\in[0,T]}$ and $(\varphi_{t\wedge T_n^j}^j(x))_{t\in[0,T]}$ are indistinguishable, we use the optional stopping theorem to obtain

$$\mathbb{E}_\mu(\exp(H_T^{T_n^j} - [H^{T_n^j}, H^{T_n^j}]_T/2)|\mathfrak{G}_{s\wedge T_n^j})$$

$$= \exp(H_{s\wedge T_n^j}^{T_n^j} - [H^{T_n^j}, H^{T_n^j}]_{s\wedge T_n^j}/2),$$

where $s \in [0, T]$. Hence,

$$\mathbb{E}_{\mu^j} \int_0^{t\wedge T_n^j} \|C(Z^j \eta)_s^{T_n^j}\|^2\, ds$$

$$= \int_0^t \mathbb{E}_\mu(\mathbf{1}_{[0,T_n^j]}(s) \eta_{s\wedge T_n^j}^2 \|CZ_{s\wedge T_n^j}^j\|^2$$

$$\times \mathbb{E}_\mu(\exp(H_T^{T_n^j} - [H^{T_n^j}, H^{T_n^j}]_T/2)|\mathfrak{G}_{s\wedge T_n^j}))\, ds.$$

Then
$$\mathbb{E}_{\mu^j} \int_0^{t \wedge T_n^j} \|C(Z^j \eta)_s^{T_n^j}\|^2 \, ds = \int_0^t \mathbb{E}_\mu(\mathbf{1}_{[0,T_n^j]}(s) \|CZ_{s \wedge T_n^j}^j\|^2) \, ds$$
$$= \mathbb{E}_\mu \int_0^{t \wedge T_n^j} \mathbf{1}_{[0,T_n^j]}(s) \|CZ_s^j\|^2 \, ds.$$

This implies $\mathbb{E}_{\mu^j} \int_0^{t \wedge T_n^j} \|C(Z^j \eta)_s^{T_n^j}\|^2 \, ds \le n$. Similarly, $\mathbb{E}_{\mu^j} \|(Z^j \eta)_s^{T_n^j}\|^2 = 1$. Then, Itô's formula shows that $\mu^j$-a.s., for all $t \in [0,T]$, $(Z^j \eta)_{t \wedge T_n^j} = \varphi^j_{t \wedge T_n^j}(x)$.

Since $\mathbb{E}_\mu \int_0^T \|CZ_t^j\|^2 \, dt < \infty$, we see that

(15) $$\mu(H_{T_n^j} - [H,H]_{T_n^j}/2 = +\infty) = 0.$$

The integration by parts formula yields
$$\mathbb{E}_\mu \left( \int_0^{T_n^j} 2\operatorname{Re}\langle Z_s^j, L_k Z_s^j \rangle \, d(Z^3)_s^k \right)^2 < \infty.$$

By (15), $\mu$ is absolutely continuous with respect to $\mu^j$. It follows that $\mu$-a.s., for all $t \in [0,T]$,

(16) $$(Z^j \eta)_{t \wedge T_n^j} = \varphi^j_{t \wedge T_n^j}(x).$$

Using again $\mathbb{E}_\mu \int_0^T \|CZ_t^j\|^2 \, dt < \infty$, we have $\mu$-a.s., for all $t \in [0,T]$, $\eta_t > 0$. Due to (16), we get $\mu$-a.s.

(17) $$\int_0^t \|C\varphi_s^j(x)/\|\varphi_s^j(x)\|\|^2 \, ds = \int_0^t \|CZ_s^j\|^2 \, ds,$$

for any $t \le T_n^j$. By $\mu(\int_0^T \|CZ_s^j\|^2 \, ds = \infty) = 0$, $t \mapsto \int_0^t \|CZ_s^j\|^2 \, ds$ is a continuous function $\mu$-a.s. Hence $\mu(T_n^j = T \text{ or } \int_0^{T_n^j} \|CZ_s^j\|^2 \, ds = n) = 1$. Applying (17) gives
$$\int_0^{T_n^j} \|C\varphi_s^j(x)/\|\varphi_s^j(x)\|\|^2 \, ds = n,$$

provided $T_n^j < T$. In addition, (17) shows that $\int_0^t \|C\varphi_s^j(x)/\|\varphi_s^j(x)\|\|^2 \, ds < n$, for any $t < T_n^j$. Therefore, $\mu$-a.s.

(18) $$T_n^j = \inf\left\{ t \in [0,T] : \int_0^t \|C\varphi_s^j(x)/\|\varphi_s^j(x)\|\|^2 \, ds \ge n \right\} \wedge T.$$

Combining the Yamada–Watanabe construction with the uniqueness of the pathwise solution of (8) given in Theorem 2, we obtain
$$\mu^1 \circ (\varphi^1(x), W^{\cdot,1})^{-1} = \mu^2 \circ (\varphi^2(x), W^{\cdot,2})^{-1}$$

on $\mathfrak{B}(C([0,T],\mathfrak{h}) \times C([0,T],\mathbb{R}^\infty))$. One way to see this is reasoning like in the proof of Theorem IX.1.7 of [53]. From (18) we get

$$
\begin{aligned}
(19) \quad (\|\varphi^1_{T^1_n}(x)\|^2 \cdot \mu^1) \circ (\varphi^1(x), W^{\cdot,1}, T^1_n)^{-1} \\
= (\|\varphi^2_{T^2_n}(x)\|^2 \cdot \mu^2) \circ (\varphi^2(x), W^{\cdot,2}, T^2_n)^{-1}
\end{aligned}
$$

on $\mathfrak{B}(C([0,T],\mathfrak{h}) \times C([0,T],\mathbb{R}^\infty) \times [0,T])$. Since

$$\|\varphi^j_{T^j_n}(x)\|^2 = (\eta_{T^j_n})^2$$
$$= \exp(-H_{T^j_n} + [H,H]_{T^j_n}/2),$$

we see that $\|\varphi^j_{T^j_n}(x)\|^2 \cdot \mu^j = \mu$. Then

$$\mu \circ (\varphi^1(x), W^{\cdot,1}, T^1_n)^{-1} = \mu \circ (\varphi^2(x), W^{\cdot,2}, T^2_n)^{-1}.$$

According to (14), (16) and (19), we have

$$\mu \circ ((Z^1)^{T^1_n}, Z^3)^{-1} = \mu \circ ((Z^2)^{T^2_n}, Z^3)^{-1}.$$

Because $\mathbb{E}_\mu(\int_0^T \|CZ^1_s\|^2 ds) < \infty$, $T^j_n \nearrow_{n\to\infty} T$ $\mu$-a.s. Letting $n \to +\infty$ yields $\mu \circ (Z^1, X^3)^{-1} = \mu \circ (Z^2, X^3)^{-1}$ on $\mathfrak{B}(C([0,T],\mathfrak{h}) \times C([0,T],\mathbb{R}^\infty) \times [0,T])$. Then, $\mathbb{Q} \circ (X,B)^{-1} = \widetilde{\mathbb{Q}} \circ (\tilde{X}, \tilde{B})^{-1}$ on $\mathfrak{B}(C([0,T],\mathfrak{h}) \times C([0,T],\mathbb{R}^\infty) \times [0,T])$. □

We now obtain easily Theorem 1 combining Propositions 1 and 2 together with Kolmogorov's theorem.

PROOF OF THEOREM 1. Proposition 2 leads to the uniqueness of the $C$-solution of (1) with initial distribution $\theta$ on $[0, +\infty[$.

By Proposition 1, there exists a solution $(\mathbb{P}_n, (X^n_t)_{t \in [0,n]}, (B^{\cdot,n}_t)_{t\in[0,n]})$ of class $C$ of (1) with initial distribution $\theta$ on $[0,n]$ for any $n \in \mathbb{N}$. Here $(X^n, B^{\cdot,n})$ takes values in $C([0,n],\mathfrak{h}) \times C([0,n],\mathbb{R}^\infty)$. Let $(X,B)$ be the coordinate mapping process defined on $\Omega = C([0,+\infty[,\mathfrak{h}) \times C([0,+\infty[,\mathbb{R}^\infty)$. Choose $\mathfrak{F}^0 = \mathfrak{B}(\Omega)$ and $\mathfrak{F}^0_t = \sigma((X_s, B_s) : s \leq t)$. For every $A \in \mathfrak{B}(C([0,n],\mathfrak{h}) \times C([0,n],\mathbb{R}^\infty))$, we define

$$\mathbb{Q}_n(\pi_n^{-1}(A)) = \mathbb{P}_n((X^n, B^n) \in A),$$

where $\pi_n : \Omega \to C([0,n],\mathfrak{h}) \times C([0,n],\mathbb{R}^\infty)$ is given by $\pi_n(\omega) = (\omega_t)_{t \in [0,n]}$. Therefore, $\mathbb{Q}_n$ is a probability measure on $(\Omega, \mathfrak{F}^0_n)$. Moreover, Proposition 2 leads to $(\mathbb{Q}_n)_{n \in \mathbb{N}}$ being a consistent family of probability measure. Since $(\Omega, \mathfrak{F}^0)$ is a standard measurable space, there exists a probability measure $\mathbb{Q}^0$ on $(\Omega, \mathfrak{F}^0)$ such that $\mathbb{Q}^0$ restricted on $\mathfrak{F}^0_n$ coincides with $\mathbb{Q}_n$ for every $n \in \mathbb{N}$. Let $(\mathfrak{F}, \mathbb{Q})$ be the completion of $(\mathfrak{F}^0, \mathbb{Q}^0)$. Let $(\mathfrak{F}_t)_{t \geq 0}$ be the usual augmentation of $(\mathfrak{F}_t)_{t \geq 0}$. Using Proposition 1 yields $(\mathbb{Q}, (X_t)_{t \geq 0}, (B_t)_{t \geq 0})$ is a $C$-solution of (1) with initial distribution $\theta$ on $[0, +\infty[$. To see this, we can use, for instance, Theorems 8.3 and 8.6 of [48]. □

**3. Regular invariant measures.** As a step toward understanding the large time behavior of the open quantum systems, this section treats Problem (ii).

3.1. *Main result.* We first apply standard arguments to obtain the Markov property of the solution of (1).

THEOREM 3. *Assume that hypotheses of Theorem 1 hold. Let $f:(\mathfrak{h}, \mathfrak{B}(\mathfrak{h})) \to (\mathbb{R}, \mathfrak{B}(\mathbb{R}))$ be a bounded measurable function. Suppose that $(\mathbb{Q}, (X_t)_{t\geq 0}, (W_t)_{t\geq 0})$ is the $C$-solution of (1). Then for any $s, t \geq 0$,*

$$\mathbb{E}(f(X_{s+t})|\mathfrak{F}_s) = \mathbb{E}(f(X_{s+t})|X_s). \tag{20}$$

*In addition, we have*

$$\mathbb{E}(f(X_{s+t})|\mathfrak{F}_s) = \int_{\mathfrak{h}} f(z) P_t(X_s, dz), \tag{21}$$

*where*

$$P_t(x, A) = \begin{cases} \mathbb{Q}_x(X_t^x \in A), & x \in Dom(C), \\ \delta_x(A), & x \notin Dom(C), \end{cases} \tag{22}$$

*provided $A \in \mathfrak{B}(\mathfrak{h})$ and $(\mathbb{Q}_x, (X_t^x)_{t\geq 0}, (B_t^{\cdot,x})_{t\geq 0})$ is the $C$-solution of (1) with initial data $x \in Dom(C)$.*

PROOF. The proof is deferred to the Appendix. □

We now adopt a set-up similar to that in Section 2 of [24], where it is discussed as the existence of quantum stationary states.

HYPOTHESIS 2. Suppose that $D$ satisfies Hypothesis 1 with $C$ substituted by $D$. Assume in addition that:

(H2.1) The set $\{x \in \mathfrak{h} : \|Dx\|^2 + \|x\|^2 \leq 1\}$ is compact in the topology of $\mathfrak{h}$.
(H2.2) There exists a vector $\hat{x}$ belonging to $Dom(D)$ such that $\|\hat{x}\| = 1$ and

$$\int_0^t \mathbb{E}\|DX_s^{\hat{x}}\|^2 \, ds \leq K_{\hat{x}} t,$$

where $K_{\hat{x}}$ is a constant depending on $\hat{x}$ and $(\mathbb{Q}, (X_t^{\hat{x}})_{t\geq 0}, (B_t^{\cdot,\hat{x}})_{t\geq 0})$ is the $D$-solution of (1) with initial data $\hat{x}$.

The next theorem provides a general criterion for the existence of a regular invariant measure for (1).

THEOREM 4. *Let $D$ satisfy Hypothesis 2. Then there exists a probability measure $\Gamma$ on the Borel $\sigma$-algebra of $\mathfrak{h}$, that is, $\mathfrak{B}(\mathfrak{h})$, such that:*

  (i) $\Gamma(Dom(D) \cap \{x \in \mathfrak{h} : \|x\| = 1\}) = 1$.
  (ii) Let $P_t(x, A)$ be given by (22) with $C$ replaced by $D$. Then

$$\Gamma(A) = \int P_t(x, A) \Gamma(dx), \tag{23}$$

for any $t \geq 0$, $x \in \mathfrak{h}$ and $A \in \mathfrak{B}(\mathfrak{h})$.
  (iii) $\int_\mathfrak{h} \|Dz\|^2 \Gamma(dz) < \infty$.

We can check directly condition (H2.2) of Hypothesis 2, for instance, in certain systems formed by an arbitrary number of identical Fermi particles (see, e.g., Section 4.2 of [47]). Nevertheless, this assumption has the disadvantage of involving the solution of (1). This motivates the following hypothesis on the operators $G$ and $L_k$ that guarantees the fulfillment of condition (H2.2).

HYPOTHESIS 3. The pair $(C, D)$ of self-adjoint positive operators in $\mathfrak{h}$ has the following properties:

(H3.1) The operator $C$ satisfies conditions (H1.2) and (H1.4) of Hypothesis 1.
(H3.2) $Dom(C) \subset Dom(D)$.
(H3.3) There exist constants $\beta \in [0, +\infty[$ and $N \in \mathbb{Z}_+$ satisfying

$$2\mathrm{Re}\langle Cx, CP_n Gx \rangle + \sum_{k=1}^\infty \|CP_n L_k x\|^2 \leq -\|Dx\|^2 + \beta(1 + \|x\|^2)$$

whenever $n \geq N$ and $x \in \mathfrak{h}_n$.
(H3.4) The operator $D$ satisfies condition (H2.1) of Hypothesis 2 and conditions (H1.1), (H1.3) and (H1.4) of Hypothesis 1 with $C$ replaced by $D$.

REMARK 6. In many situations we can find an operator $C$ such that the pair $(C, \sqrt{\alpha}C)$ obeys Hypothesis 3 for some $\alpha > 0$. This is the case, for instance, of the Jaynes–Cummings model of quantum optics and the multimode Dicke laser models (see Section 5 of [24] and Remark 2.1 of [47] for details).

REMARK 7. In Hypothesis 3 the function $x \longmapsto \|Cx\|^2$ plays the role of the Lyapunov function. In fact, the relation between the infinitesimal generator for an ordinary stochastic differential equation and the standard Lyapunov function (see, e.g., [1, 22, 38, 42]) is replaced by condition (H3.3).

The theorem below provides an intrinsic sufficient condition for the existence of regular invariant probability measures.

THEOREM 5. *Suppose that $(C, D)$ satisfies Hypothesis 3. Then there exists a probability measure $\Gamma$ on $\mathfrak{B}(\mathfrak{h})$ for which properties* (i)–(iii) *of Theorem 4 hold.*

It follows from the next lemma that condition (H2.1) of Hypothesis 2 can be expressed as the following:

(H2.1)′  $D$ has finite-dimensional spectral projections associated with bounded intervals.

Then the setting of Theorem 5 coincides essentially with the framework of Section 4 of [24], where it is treated as the existence of quantum stationary states.

LEMMA 1. *Let $D$ be a self-adjoint positive operator in $\mathfrak{h}$ such that $D$ is unbounded. Then conditions* (H2.1) *and* (H2.1)′ *are equivalent.*

REMARK 8. In Section 4, we apply Theorem 5 to a general model of harmonic oscillators with one mode. Moreover, [24] verifies essentially Hypothesis 3 in both the Jaynes–Cummings model and a multimode Dicke laser model. The analysis used in these three examples suggest to us that our main criterion for the existence of invariant measures is easy to check in a wide class of physical applications. Nevertheless, it is interesting to see how this criterion applies to other concrete quantum mechanics models.

3.2. *Proofs.* We begin by establishing that the $C$-solution of (1) depends continuously on the initial data in a distribution sense. Theorem 6 is of independent interest.

THEOREM 6. *Suppose that $C$ satisfies Hypothesis 1. Let $(x_n)_{n \in \mathbb{N}}$ be a sequence of vectors of $\mathrm{Dom}(C)$ of norm 1. Assume that there exists $x \in \mathrm{Dom}(C)$ for which $\lim_{n \to \infty} \|x_n - x\| = 0$. Let $(\mathbb{Q}_n, (X_t^{x_n})_{t \geq 0}, (B_t^{\cdot, n})_{t \geq 0})$, with $n \in \mathbb{N}$, be the $C$-solution of (1) with initial data $x_n$ (i.e., $X_0^{x_n} = x_n$ a.s.). Then for all $t \in [0, +\infty[$,*

$$(24) \qquad \mathbb{Q}_n \circ (X_t^{x_n})^{-1} \to_{n \Rightarrow \infty} \mathbb{Q} \circ (X_t^x)^{-1},$$

*provided $(\mathbb{Q}, (X_t)_{t \geq 0}, (B_t)_{t \geq 0})$ is the $C$-solution of (1) with initial data $x$. Here the symbol $\Rightarrow$ denotes weak convergence of probability measures.*

PROOF. From Propositions 1 and 2 we see that, for any $t \in [0, +\infty[$,

$$\mathbb{Q}_n \circ (X_t^{x_n})^{-1} = (\|\varphi_t(x_n)\|^2 \cdot \mathbb{P}) \circ (\varphi_t(x_n)/\|\varphi_t(x_n)\|)^{-1}$$

and $\mathbb{Q} \circ (X_t^x)^{-1} = (\|\varphi_t(x)\|^2 \cdot \mathbb{P}) \circ (\varphi_t(x)/\|\varphi_t(x)\|)^{-1}$. Here $(\varphi_t(\xi))_{t \geq 0}$ is the $C$-strong solution of (3).

Let $f : \mathfrak{h} \to \mathbb{R}$ be a bounded continuous function. Then
$$\mathbb{E}_{\mathbb{P}} f(\varphi_t(x_n)) \|\varphi_t(x_n)\|^2 - \mathbb{E}_{\mathbb{P}} f(\varphi_t(x)) \|\varphi_t(x)\|^2 = H_n^1 + H_n^2,$$
where $H_n^1 = \mathbb{E}_{\mathbb{P}} f(\varphi_t(x_n))(\|\varphi_t(x_n)\|^2 - \|\varphi_t(x)\|^2)$ and
$$H_n^2 = \mathbb{E}_{\mathbb{P}}(f(\varphi_t(x_n)) - f(\varphi_t(x))) \|\varphi_t(x)\|^2.$$
Since (3) is linear,
(25) $\quad \mathbb{E}_{\mathbb{P}} \|\varphi_t(x_n) - \varphi_t(x)\|^2 = \mathbb{E}_{\mathbb{P}} \|\varphi_t(x_n - x)\|^2 \le \|x_n - x\|^2.$

Using the Cauchy–Schwarz inequality yields
$$H_n^1 \le \sup_{z \in \mathfrak{h}} |f(z)| (\mathbb{E}_{\mathbb{P}} \|\varphi_t(x_n) - \varphi_t(x)\|^2)^{1/2} (\mathbb{E}_{\mathbb{P}}(\|\varphi_t(x_n)\| + \|\varphi_t(x)\|)^2)^{1/2}$$
$$\le \sqrt{2} \sup_{z \in \mathfrak{h}} |f(z)| (\mathbb{E}_{\mathbb{P}} \|\varphi_t(x_n) - \varphi_t(x)\|^2)^{1/2} (\|x_n\|^2 + \|(x)\|^2)^{1/2}.$$

It follows that $H_n^1 \to 0$ as $n \to \infty$.

By (25), $\mathbb{P} \circ (\varphi_t(x_n))^{-1}$ converge weakly to $\mathbb{P} \circ (\varphi_t(x))^{-1}$ as $n \to \infty$. Hence, there exists a probability space $(\widetilde{\Omega}, \widetilde{\mathfrak{F}}, \widetilde{\mathbb{P}})$ where a $\mathfrak{h}$-valued random variable $\alpha$ (resp. $\alpha_n$) is defined, with distribution $\mathbb{P} \circ (\varphi_t(x))^{-1}$ [resp. $\mathbb{P} \circ (\varphi_t(x_n))^{-1}$], and such that $\alpha_n$ converge a.s. to $\alpha$ (see, e.g., Theorem 11.7.2 of [20] or Theorem 3.1.8 of [22]). Since $\mathbb{E}_{\widetilde{\mathbb{P}}} \|\alpha\|^2 = \mathbb{E}_{\mathbb{P}} \|\varphi_t(x)\|^2 \le \|x\|^2$, Lebesgue's dominated convergence theorem leads to $H_n^2 = \mathbb{E}_{\widetilde{\mathbb{P}}}(f(\alpha_n) - f(\alpha)) \|\alpha\|^2 \longrightarrow_{n \to \infty} 0$. Therefore,

(26) $\quad (\|\varphi_t(x_n)\|^2 \cdot \mathbb{P}) \circ (\varphi_t(x_n))^{-1} \Rightarrow_{n \to \infty} (\|\varphi_t(x)\|^2 \cdot \mathbb{P}) \circ (\varphi_t(x))^{-1}.$

For any $x \in \mathfrak{h}$, set
$$\pi(x) = \begin{cases} x/\|x\|, & \text{if } x \ne 0, \\ 0, & \text{if } x = 0. \end{cases}$$
Since $(\|\varphi_t(x)\|^2 \cdot \mathbb{P}) \circ (\varphi_t(x))^{-1}(\{0\}) = 0$, (26) implies
$$((\|\varphi_t(x_n)\|^2 \cdot \mathbb{P}) \circ (\varphi_t(x_n))^{-1}) \circ \pi^{-1}$$
$$\Rightarrow_{n \to \infty} ((\|\varphi_t(x)\|^2 \cdot \mathbb{P}) \circ (\varphi_t(x))^{-1}) \circ \pi^{-1}$$
(see, e.g., Theorem 1.5.1 of [12]). This becomes
$$\int f\left(\frac{\varphi_t(x_n)}{\|\varphi_t(x_n)\|}\right) d(\|\varphi_t(x_n)\|^2 \cdot \mathbb{P})$$
$$\to_{n \to \infty} \int f\left(\frac{\varphi_t(x_n)}{\|\varphi_t(x_n)\|}\right) d(\|\varphi_t(x)\|^2 \cdot \mathbb{P})$$
whenever $f : \mathfrak{h} \to \mathbb{R}$ is a bounded continuous function, which establishes (24). □

REMARK 9. Let $\mathfrak{u} = Dom(C) \cap \{x \in \mathfrak{h} : \|x\| = 1\}$. The map from $\mathfrak{h}$ to $\mathbb{R}$ given by $x \mapsto \mathbb{E}_{\mathbb{Q}_x} f(X_t^x) \mathbf{1}_{Dom(C) \cap \mathfrak{u}}(x) + f(x)(1 - \mathbf{1}_{Dom(C) \cap \mathfrak{u}}(x))$ is measurable as soon as $f : \mathfrak{h} \to \mathbb{R}$ is bounded and continuous. This follows from Theorem 6 and Remark 3. Then, a functional form of the monotone class theorem (see, e.g., Theorem I.21 of [18]) yields that, for any $A \in \mathfrak{B}(\mathfrak{h})$, $x \mapsto P_t(x, A)$ is a measurable function from $\mathfrak{h}$ to $\mathbb{R}$.

PROOF OF THEOREM 4. Let $\hat{x}$ be as in condition (H2.2) of Hypothesis 2. Set $\mathfrak{u} = \{x \in \mathfrak{h} : \|x\| = 1\}$. For every $A \in \mathfrak{B}(\mathfrak{u})$, we define

$$\Gamma_n(A) = \frac{1}{n} \int_0^n \mathbb{Q}(X_s \in A) \, ds,$$

where $(\mathbb{Q}, (X_t)_{t \geq 0}, (B_t)_{t \geq 0})$ is the $D$-solution of (1) with initial data $\hat{x}$. Using approximations by simple functions, we can assert that

$$(27) \qquad \int_{\mathfrak{u}} g(x) \Gamma_n(dx) = \frac{1}{n} \int_0^n \mathbb{E} g(X_s) \, ds,$$

the function $g : (\mathfrak{u}, \mathfrak{B}(\mathfrak{u})) \to (\mathbb{R}, \mathfrak{B}(\mathbb{R}))$ being positive and measurable.

Set $S_l = \{x \in \mathfrak{u} : \|Dx\|^2 \leq l\}$, whenever $l > 0$. Combining Chebyshev's inequality with condition (H2.2) of Hypothesis 2, we see that

$$\Gamma_n(\mathfrak{u} \smallsetminus S_l) = \frac{1}{n} \int_0^n \mathbb{Q}(\|DX_s\|^2 > l) \, ds$$
$$\leq \frac{1}{nl} \int_0^n \mathbb{E} \|DX_s\|^2 \, ds$$
$$\leq K_{\hat{x}}/l.$$

It follows that $(\Gamma_n)_{n \in \mathbb{N}}$ is a tight family of probability measures, since $S_l$ is compact. Therefore, there exist a subsequence $(\Gamma_{n_k})_{k \in \mathbb{N}}$ and a probability measure $\Gamma$ on $\mathfrak{B}(\mathfrak{u})$ such that $\Gamma_{n_k}$ converge weakly to $\Gamma$ as $k \to \infty$ (see, e.g., Section 1.6 of [12] or Theorem 11.5.4 of [20]).

Let $f_{n,m} : \mathfrak{u} \to \mathbb{R}$ be given by $f_{n,m}(z) = \min\{\sup_{k \leq n} \|DP_k z\|^2, m\}$, with $n, m \in \mathbb{N}$. Using (27), condition (H1.4) of Hypothesis 1 with $D = C$ and condition (H2.2) of Hypothesis 2 yields

$$\int_{\mathfrak{u}} f_{n,m}(z) \Gamma_{n_k}(dz) = \frac{1}{n_k} \int_0^{n_k} \mathbb{E} f_{n,m}(X_s) \, ds$$
$$\leq \frac{1}{n_k} \int_0^{n_k} \mathbb{E} \|DX_s\|^2 \, ds$$
$$\leq K_{\hat{x}}.$$

The function $f_{n,m}$ is bounded and continuous, and so

$$\int_{\mathfrak{u}} f_{n,m}(z) \Gamma_{n_k}(dz) \longrightarrow_{k \to \infty} \int_{\mathfrak{u}} f_{n,m}(z) \Gamma(dz).$$

Hence,

$$\tag{28} \int_{\mathfrak{u}} f_{n,m}(z)\Gamma(dz) \leq K_{\hat{x}}.$$

Suppose that $\sup_{n \in \mathbb{N}} \|DP_n z\|^2 < \infty$. Then for any $y \in Dom(D)$,

$$|\langle z, Dy \rangle| = \lim_{n \to \infty} |\langle P_n z, Dy \rangle|$$
$$\leq \|y\| \sup_{n \in \mathbb{N}} \|DP_n z\|^2.$$

Hence, $z \in Dom(D^*) = Dom(D)$. By Remark 3, we have $\sup_{n \in \mathbb{N}} \|DP_n z\|^2 < \infty$ if and only if $z \in Dom(D)$. Therefore, for all $z \in \mathfrak{h}$,

$$f_{n,n}(z) \nearrow_{n \to \infty} f(z) \qquad \text{where } f(z) = \begin{cases} \|Dz\|^2, & \text{if } z \in Dom(D), \\ +\infty, & \text{if } z \notin Dom(D). \end{cases}$$

Applying the monotone convergence theorem and using (28), we get $\int_{\mathfrak{u}} f(z)\Gamma(dz) \leq K_{\hat{x}}$. Thus, $\Gamma(\mathfrak{u} \smallsetminus Dom(D)) = 0$ and $\int_{\mathfrak{u}} \|Dz\|^2 \Gamma(dz) \leq K_{\hat{x}}$.

For every probability measure $\mu$ on $\mathfrak{B}(\mathfrak{u})$, we choose

$$P_t \mu(A) = \int_{\mathfrak{u}} P_t(x, A) \mu(dx),$$

provided $A \in \mathfrak{B}(\mathfrak{u})$. Let $g : (\mathfrak{u}, \mathfrak{B}(\mathfrak{u})) \to (\mathbb{R}, \mathfrak{B}(\mathbb{R}))$ be a bounded measurable function. Approximating $g$ by simple functions gives

$$\tag{29} \int_A g(x) P_t \mu(dx) = \int_A \left( \int_{\mathfrak{u}} g(z) P_t(x, dz) \right) \mu(dx)$$

Since

$$\tag{30} \Gamma_{n_k}(\mathfrak{u} \cap Dom(D)) = \Gamma(\mathfrak{u} \cap Dom(D)) = 1,$$

the restriction of $\Gamma_{n_k}$ to $\mathfrak{B}(\mathfrak{u} \cap Dom(D))$ converge weakly to the restriction of $\Gamma$ to $\mathfrak{B}(\mathfrak{u} \cap Dom(D))$. Here $\mathfrak{u} \cap Dom(D)$ is equipped with the relative topology induced on it by $\mathfrak{h}$. Assume that $g : \mathfrak{u} \to \mathbb{R}$ is a bounded continuous function. From Theorem 6 we see that the function from $\mathfrak{u} \cap Dom(D)$ to $\mathbb{R}$ given by $x \mapsto \int_{\mathfrak{u}} g(z) P_t(x, dz)$ is continuous. Then (29) becomes

$$\int_{\mathfrak{u} \cap Dom(D)} g(z) P_t \Gamma_{n_k}(dz) = \int_{\mathfrak{u} \cap Dom(D)} \left( \int_{\mathfrak{u}} g(z) P_t(x, dz) \right) \Gamma_{n_k}(dx)$$
$$\to_{k \to \infty} \int_{\mathfrak{u} \cap Dom(D)} \left( \int_{\mathfrak{u}} g(z) P_t(x, dz) \right) \Gamma(dx)$$
$$= \int_{\mathfrak{u} \cap Dom(D)} g(z) P_t \Gamma(dz).$$

From (30), we obtain that (see, e.g., Corollary 3.3.2 of [22])

$$\tag{31} P_t \Gamma_{n_k} \Rightarrow_{k \to \infty} P_t \Gamma \qquad \text{on } \mathfrak{u}.$$

Using (21) and (27) yields

$$P_t \Gamma_{n_k}(A) = \frac{1}{n_k} \int_0^{n_k} \mathbb{E} P_t(X_s, A) \, ds$$

$$= \frac{1}{n_k} \int_0^{n_k} \mathbb{E}(\mathbb{E}(1_A(X_{s+t})|\mathfrak{F}_s)) \, ds,$$

where $A \in \mathfrak{B}(\mathfrak{u})$. It follows that $P_t \Gamma_{n_k}(A) = \frac{1}{n_k} \int_0^{n_k} \mathbb{P}(X_{s+t} \in A) \, ds$. Hence, for any close set $F$ in $\mathfrak{u}$,

$$\limsup_{k \to \infty} P_t \Gamma_{n_k}(F)$$
$$= \limsup_{k \to \infty} \left( \Gamma_{n_k}(F) + \frac{1}{n_k} \left( \int_{n_k}^{n_k+t} \mathbb{P}(X_s \in F) \, ds - \int_0^t \mathbb{P}(X_s \in F) \, ds \right) \right)$$
$$= \limsup_{k \to \infty} \Gamma_{n_k}(F).$$

Applying the Portmanteau theorem, we see that

$$\limsup_{k \to \infty} P_t \Gamma_{n_k}(F) \leq \Gamma(F).$$

Therefore, $P_t \Gamma_{n_k} \Rightarrow_{k \to \infty} \Gamma$ on $\mathfrak{u}$. From (31) we have $P_t \Gamma = \Gamma$. For any $A \in \mathfrak{B}(\mathfrak{h})$, we set $\Gamma(A) = \Gamma(A \cap \mathfrak{u})$, and the theorem follows. $\square$

To prove Theorem 5, we need the following lemma which provides global estimates for $\mathbb{E}\|DX_t\|^2$.

LEMMA 2. *Let the hypotheses of Theorem 5 hold. Suppose that $\theta$ is a probability measure on $\mathfrak{B}(\mathfrak{h})$ such that $\theta(Dom(C) \cap \{x \in \mathfrak{h} : \|x\| = 1\}) = 1$ and $\int_{\mathfrak{h}} \|Cx\|^2 \theta(dx) < \infty$. Let $(\mathbb{Q}, (X_t)_{t \geq 0}, (B_t)_{t \geq 0})$ be the $C$-solution of (1) with initial distribution $\theta$. Then for all $t \geq 0$,*

$$(32) \qquad \int_0^t \mathbb{E}_{\mathbb{Q}} \|DX_s\|^2 \, ds \leq \mathbb{E}_{\mathbb{Q}} \|CX_0\|^2 + 2\beta t.$$

PROOF. Assume that $\xi$ is distributed according to $\theta$. Let $\psi_{\cdot,n}(\xi)$ be the continuous strong solution of the following stochastic differential equation on $\mathfrak{h}_n$:

$$(33) \quad \psi_{t,n}(\xi) = P_n \xi + \int_0^t P_n G \psi_{s,n}(\xi) \, ds + \sum_{k=1}^n \int_0^t P_n L_k \psi_{s,n}(\xi) \, dW_s^k.$$

Applying Itô's formula yields (see, e.g., [47])

$$(34) \qquad \mathbb{E}\|\psi_{t,n}(\xi)\|^2 \leq \mathbb{E}\|\xi\|^2.$$

Combining Itô's formula with condition (H3.3) of Hypothesis 3, we obtain that, for any $n \geq N$,

$$\mathbb{E}\|C\psi_{t,n}(\xi)\|^2$$
$$\leq \mathbb{E}\|C\psi_{0,n}(\xi)\|^2 + \mathbb{E}\int_0^t (-\|D\psi_{s,n}(\xi)\|^2 + \beta(1 + \|\psi_{t,n}(\xi)\|^2))\, ds.$$

Then, (34) becomes

$$(35) \qquad \int_0^t \mathbb{E}\|D\psi_{s,n}(\xi)\|^2 \, ds \leq \mathbb{E}\|CP_n\xi\|^2 + 2\beta T$$

for all $n \geq N$.

Let $\nu$ be the Lebesgue measure on $\mathfrak{B}([0,T])$, where $T$ is a positive real number. Because the unit ball of $L^2(\Omega \times [0,T], \mathbb{P} \otimes \nu; (Dom(D), \langle \cdot, \cdot \rangle_D))$ is weak* compact, (34) and (35) imply that there exists a subsequence $(\psi_{\cdot,n_l}(\xi))_{l \in \mathbb{N}}$ such that $\psi_{\cdot,n_l}(\xi)$ converge weakly in $L^2(\Omega \times [0,T], \mathbb{P} \otimes \nu; (Dom(D), \langle \cdot, \cdot \rangle_D))$ as $l \to \infty$. From Section 3.2 of [47] we have $\mathbb{E}\|\psi_{t,n_l}(\xi) - \varphi_t(\xi)\|^2 \to_{l \to \infty} 0$, and so $\psi_{\cdot,n_l}(\xi)$ converge weakly to $\varphi_\cdot(\xi)$ as $l \to \infty$ in $L^2(\Omega \times [0,T], \mathbb{P} \otimes \nu; (Dom(D), \langle \cdot, \cdot \rangle_D))$. Hence,

$$\int_0^T \mathbb{E}\|\varphi_t(\xi)\|_D^2 \, dt \leq \liminf_{l \to \infty} \int_0^T \mathbb{E}\|\psi_{t,n_l}(\xi)\|_D^2 \, dt.$$

By (34) and property (i) of Theorem 2,

$$(36) \qquad \int_0^T \mathbb{E}\|D\varphi_t(\xi)\|^2 \, dt \leq \liminf_{l \to \infty} \int_0^T \mathbb{E}\|D\psi_{t,n_l}(\xi)\|^2 \, dt.$$

According to (35) and (36), we have $\int_0^T \mathbb{E}\|D\varphi_t(\xi)\|^2 \, dt \leq \mathbb{E}\|C\xi\|^2 + 2\beta T$. Applying Propositions 1 and 2, we get (32). □

PROOF OF THEOREM 5.  Since $Dom(C) \subset Dom(D)$, the $C$-solution of (1) with initial distribution $\theta$ coincides with the $D$-solution of (1) with initial distribution $\theta$. Combining Theorem 4 with Lemma 2, we obtain the assertion of the theorem. □

PROOF OF LEMMA 1.  Suppose that assumption (H2.1)′ holds. We define $E$ to be the spectral decomposition of $D$ (see, e.g., [55]). Then the following property holds.

PROPERTY P.  Every bounded interval $I$, with $E(I) \neq 0$, contains a point $p$ such that

$$(37) \qquad\qquad\qquad E(\{p\}) \neq 0. \qquad\qquad\qquad\qquad \square$$

In fact, there exists a nonzero vector $y$ in the range of $E(I)$ since $E(I) \neq 0$. We write $E_{y,y}$ for the finite positive measure $\langle y, Ey \rangle$. Take $I_1 = I$. We define the interval $I_{n+1}$, with $n \in \mathbb{N}$, by the recurrence relation $I_{n+1} = I_n^{j_n}$, provided $j_n \in \{0,1\}$ satisfies

(38) $$E_{y,y}(I_n^{j_n}) \geq E_{y,y}(I_n^{1-j_n}).$$

Here $I_n^0$, $I_n^1$ is a partition of $I_n$ into two disjoint subintervals of the same length.

Condition (H2.1)' implies that the dimension of the range of $E(I)$ is equal to a natural number $N$. Then $\mathcal{I} = \{n \in \mathbb{N} : E_{y,y}(I_n) \neq E_{y,y}(I_{n+1})\}$ has at most $N-1$ elements. Conversely, suppose that $n_1, \ldots, n_N$ belong to $\mathcal{I}$. By

(39) $$E_{y,y}(I_n) = E_{y,y}(I_n^0) + E_{y,y}(I_n^1),$$

for any $n \in \mathbb{N}$, we have $E_{y,y}(I_n^{1-j_n}) \neq 0$ whenever $n \in \{n_1, \ldots, n_N\}$. Since $I_{n_1}^{1-j_{n_1}}, \ldots, I_{n_N}^{1-j_{n_N}}, I_{n_N}^{j_{n_N}}$ are disjoint, $E(I_{n_1}^{1-j_{n_1}}), \ldots, E(I_{n_N}^{1-j_{n_N}})$ and $E(I_{n_N}^{j_{n_N}})$ are orthogonal to each other. It follows that the dimension of the range of $E(I)$ is greater than or equal to $N+1$, which is impossible. Therefore, there exists $N_0 \in \mathbb{N}$ such that $E_{y,y}(I_n) = E_{y,y}(I_{N_0})$ for all $n \geq N_0$.

Applying Lebesgue's dominated convergence theorem yields

$$E_{y,y}\left(\bigcap_{n \in \mathbb{N}} I_n\right) = \lim_{n \to \infty} \int_{\mathbb{R}} \mathbf{1}_{I_n}(t) E_{y,y}(dt)$$
$$= E_{y,y}(I_{N_0}).$$

Combining (38) with (39), we obtain $E_{y,y}(I_{N_0}) > 0$. Thus, $\bigcap_{n \in \mathbb{N}} I_n = \{p\}$. From

$$\|E(\{p\})y\|^2 = \langle y, E(\{p\})y \rangle$$
$$= E_{y,y}(\{p\}),$$

we have (37).

Combining the definition of the resolution of the identity (see, e.g., Definition 12.17 of [55]) with repeated application of Property P, we see that there exist real numbers $p_j$, where $j$ are natural numbers less than $M \in \mathbb{N} \cup \{+\infty\}$, such that, for any $k \in \mathbb{N}$, the number of elements of $\{j : p_j \in [0,k]\}$ is at most the dimension of the range of $E([0,k])$ and $E([0,k]) = \sum_{p_j \in [0,k]} E(\{p_j\})$. By the spectrum of $D$ is concentrated on $[0, +\infty[$, for any $x \in Dom(D)$ and $y \in \mathfrak{h}$, we have

$$\langle y, Dx \rangle = \sum_{j < M} p_j \langle y, E(\{p_j\})x \rangle.$$

Using the Banach–Steinhaus theorem yields $Dx = \sum_{j < M} p_j E(\{p_j\})x$ (see, e.g., Theorem 12.6 of [55]). The operator $D$ is unbounded, and so $M = \infty$. Then $p_j \to \infty$ as $j \to \infty$, which yields $D$ satisfies condition (H2.1) of Hypothesis 2.

On the other hand, assume that condition (H2.1) holds. Then $(D^2 + I)^{-1}$ is a compact self-adjoint operator in $\mathfrak{h}$. Since 0 is not an eigenvalue of $(D^2 + I)^{-1}$, there exist an orthonormal basis $(e_n)_{n \in \mathbb{Z}_+}$ of $\mathfrak{h}$ and a sequence of positive real numbers $(\beta_n)_{n \in \mathbb{Z}_+}$ for which $\lim_{n \to \infty} \beta_n = 0$ and $(D^2 + I)^{-1} = \sum_{n \in \mathbb{Z}_+} \beta_n \langle e_n, \cdot \rangle e_n$ (see, e.g., Theorem 19B of [62]). Hence, $\beta_n \in ]0,1]$, provided $n \in \mathbb{Z}_+$, and

$$D = \sum_{n \in \mathbb{Z}_+} \sqrt{1/\beta_n - 1} \langle e_n, \cdot \rangle e_n.$$

This implies assumption (H2.1)'. □

**4. Application.** This section illustrates the main results of this paper with a one-dimensional harmonic oscillator.

THEOREM 7. *Let the assumptions of Example 1 hold. Suppose that $p$ is a natural number greater than or equal to 4. Assume that $\theta$ is a probability measure on $\mathfrak{B}(\mathfrak{h})$ such that $\theta(Dom(N^p) \cap \{x \in \mathfrak{h} : \|x\| = 1\}) = 1$ and $\int_{\mathfrak{h}} \|N^p x\|^2 \theta(dx) < \infty$. Then, (1) has a unique solution of class $N^p$ with initial distribution $\theta$ provided $|\alpha_4| \geq |\alpha_5|$.*

PROOF. Let $(e_j)_{j \in \mathbb{Z}_+}$ be the canonical orthonormal basis on $l^2(\mathbb{Z}_+)$. Then $Dom(N^p) = \{x \in \mathfrak{h} : \sum_{j=0}^{\infty} j^{2p} |\langle e_j, x \rangle|^2 < \infty\}$ and $N^p = \sum_{j=0}^{\infty} j^p \langle e_j, \cdot \rangle$. It follows that $N^p$ is a self-adjoint positive operator in $l^2(\mathbb{Z}_+)$ and the conditions (H1.2) and (H1.4) hold. Furthermore, the set $\{x \in l^2(\mathbb{Z}_+) : \|N^p x\|^2 + \|x\|^2 \leq 1\}$ is compact since $\lim_{j \to \infty} j^p = \infty$.

By $Dom(N^4) \subset Dom(H) \cap \bigcap_{k=1}^{6} (Dom(L_k^*) \cap Dom(L_k^* L_k))$, the operator $N^p$ satisfies hypothesis (H1.1).

For simplicity, assume that $n \geq 2$. A long easy computation shows that, for any $x$ belonging to the linear span of $e_0, \ldots, e_n$,

$$(40) \qquad 2\mathrm{Re}\langle N^p x, N^p P_n G x\rangle + \sum_{k=1}^{6} \|N^p P_n L_k x\|^2 \leq \sum_{j=0}^{n} c_j |x_j|^2,$$

where $x_j = \langle e_j, x \rangle$ and

$$\begin{aligned}
c_j &= |\beta_1|(2p(\sqrt{j} + \sqrt{j+1})j^{2p-1} + \sqrt{j} P_{2p-2}(j) + \sqrt{j+1} P_{2p-2}(j)) \\
&\quad + |\alpha_1|^2(-2pj^{2p} + p(2p-1)j^{2p-1} + P_{2p-2}(j)) \\
&\quad + |\alpha_2|^2(2pj^{2p} + p(2p+1)j^{2p-1} + P_{2p-2}(j)) \\
&\quad + |\alpha_4|^2(-4pj^{2p+1} + 8p^2 j^{2p} + P_{2p-1}(j)) \\
&\quad + |\alpha_5|^2(4pj^{2p+1} + 8p(p+1)j^{2p} + P_{2p-1}(j)).
\end{aligned}$$

Here we use the same symbol $P_l$ for different polynomials of degree $l$ whose coefficients depend only on $l$. It follows that condition (H1.3) holds when $|\alpha_4| \geq |\alpha_5|$. Theorem 1 now leads to our claim. □

THEOREM 8. *Let the assumptions of Example 1 hold. Assume that $p$ is a natural number greater than or equal to 4. Suppose that either $|\alpha_4| > |\alpha_5|$ or $|\alpha_4| = |\alpha_5|$ with $|\alpha_2|^2 - |\alpha_1|^2 + 4(2p+1)|\alpha_4|^2 < 0$. Then there exists a probability measure $\Gamma$ on $\mathfrak{B}(\mathfrak{h})$ satisfying (23) such that $\int_{\mathfrak{h}} \|N^p\|^2 \Gamma(dz) < \infty$ and*

$$\Gamma(Dom(N^p) \cap \{x \in \mathfrak{h} : \|x\| = 1\}) = 1.$$

PROOF. We return to the proof of Theorem 7. Let $|\alpha_4| > |\alpha_5|$. By (40),

$$2\mathrm{Re}\langle N^p x, N^p P_n G x\rangle + \sum_{k=1}^{6} \|N^p P_n L_k x\|^2$$

$$\leq K|x_0|^2 + \sum_{j=1}^{n} |x_j|^2 j^{2p}(4pj(|\alpha_5|^2 - |\alpha_4|^2) + \mathcal{O}_{2p}(j)/j^{2p}),$$

where $K$ is a positive constant and $n \geq 2$. Here $(O_{2p}(j))_{j \in \mathbb{Z}_+}$ is a sequence for which $\lim_{j\to\infty} \mathcal{O}_{2p}(j)/j^{2p}$ exists. For any $r > 0$,

$$\lim_{j \to \infty} (rj(|\alpha_5|^2 - |\alpha_4|^2) + O_{2p}(j)/j^{2p}) = -\infty,$$

and so the pair $(C, \sqrt{\alpha}C)$ obeys condition (H2.3) for any $\alpha \in\, ]0, 4p(|\alpha_4|^2 - |\alpha_5|^2)[$.

Suppose that $|\alpha_4| = |\alpha_5|$ and $|\alpha_2|^2 - |\alpha_1|^2 + 4(2p+1)|\alpha_4|^2 < 0$. From (40), we obtain

$$2\mathrm{Re}\langle N^p x, N^p P_n G x\rangle + \sum_{k=1}^{6} \|N^p P_n L_k x\|^2$$

$$\leq K|x_0|^2 + \sum_{j=1}^{n} |x_j|^2 j^{2p}(8p(2p+1)|\alpha_4|^2 + 2p(|\alpha_2|^2 - |\alpha_1|^2)$$

$$+ o_{2p}(j)/j^{2p}),$$

where $n \geq 2$ and $K > 0$. Here $(o_{2p}(j))_{j \in \mathbb{Z}_+}$ is a sequence such that

$$\lim_{j \to \infty} o_{2p}(j)/j^{2p} = 0.$$

Therefore, $(C, \sqrt{\alpha}C)$ satisfies the condition (H2.3) for any

$$\alpha \in\, ]0, 2p(|\alpha_1|^2 - |\alpha_2|^2 - 4(2p+1)|\alpha_4|^2)[.$$

From Theorem 5 we now obtain the claim of this theorem. □

A particular case of Theorem 8 applies in the following simple damped harmonic oscillator.

EXAMPLE 3. In the setting of Example 1, consider $\beta_1 = 0$, $\beta_2 = \omega$, $\beta_3 = 0$, $\alpha_1 = \sqrt{A(\nu+1)}$, $\alpha_2 = \sqrt{A\nu}$ and $\alpha_j = 0$, $(3 \leq j \leq 6)$. Here $\omega$, $A$ and $\nu$ are positive real numbers.

Example 3 describes a mode of the quantized radiation field of an ideal resonator which interacts with two-level atoms that pass through the resonator (see, e.g., [21]). In this situation, $A$ is the energy decay rate, $\nu$ is the number of thermal excitations in the steady state and $\omega$ is the natural (circular) frequency. Since $|\alpha_2|^2 - |\alpha_1|^2 = -A < 0$, Theorem 8 gives the existence of a regular stationary measure. This is in agreement with Section 4.1 of [52], where the existence of a unique faithful stationary state is studied.

Theorem 8 also covers the next basic radiation-matter interaction mechanics.

EXAMPLE 4. In the context of Example 1, we define $\beta_1$, $\beta_2$, $\alpha_1$, $\alpha_2$, $\alpha_3$ and $\alpha_6$ to be equal to 0. Moreover, we set $\beta_3 \in \mathbb{R}$, $\alpha_4 > 0$ and $\alpha_5 \geq 0$.

Example 4 simulates a two-photon absorption and emission process. Since the phenomenon of two-photon absorption was observed by Kaiser and Garret in [37], models like Example 4 have been discussed in the physical literature (see, e.g., [54]). Using Theorem 8 yields the existence of a regular invariant measure whenever $\alpha_4 > \alpha_5$, which is in agreement with [27]. In [27], Fagnola and Quezada characterized all the invariant states corresponding to Example 4 with $\alpha_4 > \alpha_5$.

## APPENDIX

We can prove the Markov property of the $C$-solution of (1) using techniques of well-posed martingale problems.

PROOF OF THEOREM 3. To prove (20), we modify the proof of Theorem 4.4.2(a) of [22]. Let $s \geq 0$. Consider the set $A$ in $\mathfrak{F}_s$ such that $\mathbb{Q}(A) > 0$. For any $B \in \mathfrak{F}$, we define $\mathbb{Q}_1(B) = \mathbb{E}(1_A \mathbb{E}(1_B|\mathfrak{F}_s))/\mathbb{Q}(A)$ and $\mathbb{Q}_2(B) = \mathbb{E}(1_A \mathbb{E}(1_B|X_s))/\mathbb{Q}(A)$. Then $\mathbb{Q}_1, \mathbb{Q}_2 \ll \mathbb{Q}$, that is, $\mathbb{Q}_1, \mathbb{Q}_2$ are absolutely continuous with respect to $\mathbb{Q}$.

Let $r \geq 0$. For any bounded measurable function $g:(\Omega, \mathfrak{F}_{s+r}) \to (\mathbb{R}, \mathfrak{B}(\mathbb{R}))$, we have $\mathbb{E}((W^k_{s+t} - W^k_{s+r})g|\mathfrak{F}_{s+r}) = 0$, provided $t \geq r$ and $k \in \mathbb{N}$. Hence, for any $j = 1, 2$ and $k \in \mathbb{N}$, $\mathbb{E}_{\mathbb{Q}_j}((W^k_{s+t} - W^k_{s+r})g) = 0$, and so

$$(41) \qquad \mathbb{E}_{\mathbb{Q}_j}(W^k_{s+t}|\mathfrak{F}_{s+r}) = W^k_{s+r}.$$

Write $B_t^k = W_{s+t}^k - W_t^k$ whenever $t \geq 0$ and $k \in \mathbb{N}$. Since $\mathbb{Q}_j \ll \mathbb{Q}$, $\mathbb{Q}_j$-a.s., for all $t \in [0, +\infty[$,

$$X_{s+t} = X_s + \int_s^{s+t} G(X_r)\,dr + \sum_{k=1}^\infty \int_s^{s+t} L_k(X_r)\,dB_r^k.$$

Using (41), $\mathbb{Q}_j \ll \mathbb{Q}$ and the Lévy characterization of Brownian motion, we obtain that $B^1, B^2, \ldots,$ is a sequence of independent Brownian motions on $(\Omega, \mathfrak{F}, (\mathfrak{F}_{s+t})_{t\geq 0}, \mathbb{Q}_j)$. Therefore, $(\Omega, \mathfrak{F}, (\mathfrak{F}_{s+t})_{t\geq 0}, \mathbb{Q}_j, (X_{s+t})_{t\geq 0}, (B_t^k)_{t\geq 0}^{k\in\mathbb{N}})$ is a $C$-solution of (1) with initial data distributed according to the law of $X_s$. By

$$\mathbb{Q}_1 \circ (X_s)^{-1} = \mathbb{Q}_2 \circ (X_s)^{-1} = \mathbb{Q}(X_s \in \cdot | A),$$

Theorem 1 leads to $\mathbb{E}_{\mathbb{Q}_1}(f(X_{s+t})) = \mathbb{E}_{\mathbb{Q}_2}(f(X_{s+t}))$ for any $t \geq 0$. This gives (20).

Let $Q$ be the regular conditional distribution for $(X, W)$ with range space $C([0, +\infty[, \mathfrak{h}) \times C([0, +\infty[, \mathbb{R}^\infty)$ given $(X_0, W_0)$. From Theorem 1 we deduce that $\mathbb{Q}$-a.s. for all $\omega \in \Omega$, $\mathbb{Q}_{X_0(\omega)} \circ (X^{X_0(\omega)}, B^{X_0(\omega)})^{-1} = Q(\omega, \cdot)$. This follows paraphrasing the proof of Proposition IX.1.4 of [53]. Then for $\omega$ in a set of probability 1 for $\mathbb{Q}$, $P_t(X_0(\omega), A) = Q(\omega, \pi^{-1}(A))$, where $A \in \mathfrak{B}(\mathfrak{h})$ and the map $\pi$ is defined by $\pi(a, b) = a(t)$ for any $a \in C([0, +\infty[, \mathfrak{h})$ and $b \in C([0, +\infty[, \mathbb{R}^\infty)$. Using the definition of regular conditional distribution yields

(42) $$\mathbb{E}(1_A(X_t)|X_0) = P_t(X_0, A), \qquad \mathbb{Q}\text{-a.s.}$$

for any $A \in \mathfrak{B}(\mathfrak{h})$.

Suppose that $s$ is greater than 0. Set $B_t^k = W_{s+t}^k - W_s^k$. Then $(B_t^k)_{t\geq 0}^{k\in\mathbb{N}}$ is a sequence of independent $(\mathfrak{F}_{s+t})_{t\geq 0}$-Brownian motions. Furthermore, we have $(\Omega, \mathfrak{F}, (\mathfrak{F}_{s+t})_{t\geq 0}, \mathbb{Q}, (X_{s+t})_{t\geq 0}, (B_t^k)_{t\geq 0}^{k\in\mathbb{N}})$ is the $C$-solution of (1) with initial data distributed according to the law of $X_s$. According to (42),

$$\mathbb{E}(1_A(X_{s+t})|X_s) = P_t(X_s, A), \qquad \mathbb{Q}\text{-a.s.}$$

whenever $A \in \mathfrak{B}(\mathfrak{h})$, and so (21) follows. $\square$

**Acknowledgments.** The authors wishes to express their thanks to Franco Fagnola and Alberto Barchielli for helpful comments.

## REFERENCES

[1] Arnold, L. (1974). *Stochastic Differential Equations*: *Theory and Applications*. Wiley, New York. MR0443083

[2] Accardi, L., Fagnola, F. and Hachicha, S. (2006). Generic $q$-Markov semigroups and speed of convergence of q-algorithms. *Infin. Dimens. Anal. Quantum Probab. Relat. Top.* **9** 567–594. MR2282720

[3] ALICKI, R. and FANNES, M. (1987). Dilations of quantum dynamical semigroups with classical Brownian motion. *Comm. Math. Phys.* **108** 353–361. MR0874898

[4] BARCHIELLI, A. and BELAVKIN, V. P. (1991). Measurements continuous in time and a posteriori states in quantum mechanics. *J. Phys. A Math. Gen.* **24** 1495–1514. MR1121822

[5] BARCHIELLI, A. and HOLEVO, A. S. (1995). Constructing quantum measurement processes via classical stochastic calculus. *Stochastic Process. Appl.* **58** 293–317. MR1348380

[6] BARCHIELLI, A., PAGANONI, A. M. and ZUCCA, F. (1998). On stochastic differential equations and semigroups of probability operators in quantum probability. *Stochastic Process. Appl.* **73** 69–86. MR1603834

[7] BARCHIELLI, A. and PAGANONI, A. M. (2003). On the asymptotic behaviour of some stochastic differential equations for quantum states. *Infin. Dimens. Anal. Quantum Probab. Relat. Top.* **6** 223–243. MR1991493

[8] BELAVKIN, V. P. (1989). Nondemolition measurements, nonlinear filtering and dynamic programming of quantum stochastic processes. *Lecture Notes in Control and Inform. Sci.* **121** 245–265. Springer, Berlin. MR1231124

[9] BELAVKIN, V. P. and KOLOKOLTSOV, V. N. (1991). Quasiclassical asymptotics of quantum stochastic equations. *Teoret. Mat. Fiz.* **89** 163–177. [Translation in *Theoret. Math. Phys.* **89** (1992) 1127–1138.] MR1151379

[10] BELAVKIN, V. P. (1992). Quantum stochastic calculus and quantum nonlinear filtering. *J. Multivar. Anal.* **42** 171–201. MR1183841

[11] BELAVKIN, V. P. (1993, 1994). Quantum diffusion, measurement and filtering. I, II. *Theor. Probab. Appl.* **38**, **39** 573–585, 363–378. MR1317995, MR1347181

[12] BILLINGSLEY, P. (1968). *Convergence of Probability Measures.* Wiley, New York. MR0233396

[13] CARMICHAEL, H. (1993). *An Open Systems Approach to Quantum Optics.* Springer, Berlin.

[14] CHEBOTAREV, A. M. and FAGNOLA, F. (1998). Sufficient conditions for conservativity of minimal quantum dynamical semigroups. *J. Funct. Anal.* **153** 382–404. MR1614586

[15] CHEBOTAREV, A. M. (2000). *Lectures on Quantum Probability.* Sociedad Matemática Mexicana, Mexico. MR1925129

[16] COHEN-TANNOUDJI, C., DIU, B. and LALOË, F. (1977). *Quantum Mechanics* **1**. Hermann, Paris.

[17] DALIBARD, J. and CASTIN, Y. (1992). Wave-function approach to dissipative processes in quantum optics. *Phys. Rev. Lett.* **68** 580–583.

[18] DELLACHERIE, C. and MEYER, P.-A. (1978). *Probabilities and Potential.* North-Holland, Amsterdam. MR0521810

[19] DIÓSI, L. (1988). Continuous quantum measurement and Itô formalism. *Phys. Lett. A* **129** 419–423. MR0952519

[20] DUDLEY, R. M. (1989). *Real Analysis and Probability.* Wadsworth, Pacific Grove, CA. MR0982264

[21] ENGLERT, B.-G. and MORIGI, G. (2002). Five lectures on dissipative master equations. In *Coherent Evolution in Noisy Environments. Lecture Notes in Phys.* **611** 55–106. Springer, Berlin. MR2012175

[22] ETHIER, S. N. and KURTZ, T. G. (1986). *Markov Processes. Characterization and Convergence.* Wiley, New York. MR0838085

[23] FAGNOLA, F. (1999). Quantum Markov semigroups and quantum flows. *Proyecciones* **18** 1–144. MR1814506

[24] FAGNOLA, F. and REBOLLEDO, R. (2001). On the existence of stationary states for quantum dynamical semigroups. *J. Math. Phys.* **42** 1296–1308. MR1814690
[25] FAGNOLA, F. and REBOLLEDO, R. (2002). Subharmonic projections for a quantum Markov semigroup. *J. Math. Phys.* **43** 1074–1082. MR1878987
[26] FAGNOLA, F. and REBOLLEDO, R. (2003). Transience and recurrence of quantum Markov semigroups. *Probab. Theory Related Fields* **126** 289–306. MR1990058
[27] FAGNOLA, F. and QUEZADA, R. (2005). Two-photon absorption and emission process. *Infin. Dimens. Anal. Quantum Probab. Relat. Top.* **8** 573–591. MR2184084
[28] GARDINER, C. W. and ZOLLER, P. (2000). *Quantum Noise.* Springer, Berlin. MR1736115
[29] GATAREK, D. and GISIN, N. (1991). Continuous quantum jumps and infinite-dimensional stochastic equations. *J. Math. Phys.* **32** 2152–2157. MR1123608
[30] GHOSE, S., ALSING, P., SANDERS, B. C. and DEUTSCH, I. (2005). Entanglement and the quantum-to-classical transition. *Phys. Rev. A* **72** 014102. MR2168054
[31] GISIN, N. and PERCIVAL, I. C. (1992). The quantum-state diffusion model applied to open systems. *J. Phys. A Math. Gen.* **25** 5677–5691. MR1192024
[32] GOUGH, J. and SOBOLEV, A. (2004). Continuous measurement of canonical observables and limit stochastic Schrodinger equations. *Phys. Rev. A* **69** 032107. MR2062685
[33] HUDSON, R. L. and PARTHASARATHY, K. R. (1984). Quantum Itô's formula and stochastic evolutions. *Comm. Math. Phys.* **93** 301–323. MR0745686
[34] HOLEVO, A. S. (1996). On dissipative stochastic equations in a Hilbert space. *Probab. Theory Related Fields* **104** 483–500. MR1384042
[35] HOLEVO, A. S. (1996). Stochastic differential equations in Hilbert space and quantum Markovian evolutions. In *Probability Theory and Mathematical Statistics* (*Tokyo, 1995*) 122–131. World Science, River Edge, NJ. MR1467932
[36] IKEDA, N. and WATANABE, S. (1981). *Stochastic Differential Equations and Diffusion Processes.* North-Holland, Amsterdam. MR1011252
[37] KAISER, W. and GARRET, C. G. B. (1961). Two-photon excitation in $CaF_2:Eu^{2+}$. *Phys. Rev. Lett.* **7** 229–231.
[38] KHAS'MINSKII, R. Z. (1980). *Stochastic Stability of Differential Equations.* Sijthoff and Nordhoff, Alphen aan den Rijn. MR0600653
[39] KOLOKOLTSOV, V. N. (1998). Long time behavior of continuously observed and controlled quantum systems (a study of the Belavkin quantum filtering equation). *Quantum Probab. Comm.* **10** 229–243. MR1689486
[40] KOLOKOLTSOV, V. N. (1998). Localization and analytic properties of the solutions of the simplest quantum filtering equation. *Rev. Math. Phys.* **10** 801–828. MR1643961
[41] KOLOKOLTSOV, V. N. (2000). *Semiclassical Analysis for Diffusions and Stochastic Processes.* Springer, Berlin. MR1755149
[42] MATTINGLY, J. C., STUART, A. M. and HIGHAM, D. J. (2002). Ergodicity for SDEs and approximations: Locally Lipschitz vector fields and degenerate noise. *Stochastic Process. Appl.* **101** 185–232. MR1931266
[43] MEYER, P. A. (1993). *Quantum Probability for Probabilists.* Springer, Berlin. MR1222649
[44] MIKULEVICIUS, R. and ROZOVSKII, B. L. (1999). Martingale problems for stochastic PDE's. In *Stochastic Partial Differential Equations*: *Six Perspectives* 243–325. Amer. Math. Soc., Providence, RI. MR1661767
[45] MORA, C. M. (2004). Numerical simulation of stochastic evolution equations associated to quantum Markov semigroups. *Math. Comp.* **73** 1393–1415. MR2047093

[46] MORA, C. M. (2005). Numerical solution of conservative finite-dimensional stochastic Schrödinger equations. *Ann. Appl. Probab.* **15** 2144–2171. MR2152256
[47] MORA, C. M. and REBOLLEDO, R. (2007). Regularity of solutions to linear stochastic Schrödinger equations. *Infin. Dimens. Anal. Quantum Probab. Relat. Top.* **10** 237–259. MR2337521
[48] ONDREJÁT, M. (2004). Uniqueness for stochastic evolution equations in Banach spaces. *Dissertationes Math.* **426** 63. MR2067962
[49] PARTHASARATHY, K. R. (1992). *An Introduction to Quantum Stochastic Calculus*. Birkhäuser, Basel. MR1164866
[50] PROTTER, P. (2004). *Stochastic Integration and Differential Equations*, 2nd ed. Springer, Berlin. MR2020294
[51] PERCIVAL, I. C. (1998). *Quantum State Diffusion*. Cambridge Univ. Press. MR1666826
[52] REBOLLEDO, R. (2005). A view on decoherence via master equations. *Open Syst. Inf. Dyn.* **12** 37–54. MR2135004
[53] REVUZ, D. and YORK, M. (1999). *Continuous Martingales and Brownian Motion*, 3rd ed. Springer, Berlin. MR1725357
[54] ROA, L. (1994). Phase squeezing in a three-level atom micromaser. *Phys. Rev. A* **50** R1995–R1998.
[55] RUDIN, W. (1973). *Functional Analysis*. McGraw-Hill, New York. MR0365062
[56] SAKURAI, J. J. (1994). *Modern Quantum Mechanics*, rev. ed. Addison-Wesley, Reading, MA.
[57] SCOTT, A. J. and MILBURN, G. J. (2001). Quantum nonlinear dynamics of continuously measured systems. *Phys. Rev. A* **63** 42101. MR1812792
[58] SCHACK, R., BRUN, T. A. and PERCIVAL, I. C. (1995). Quantum state diffusion, localization and computation. *J. Phys. A* **28** 5401–5413. MR1364146
[59] SPEHNER, D. and ORSZAG, M. (2002). Temperature enhanced squeezing in cavity QED. *J. Opt. B* **4** 326–335.
[60] UMANITÀ, V. (2006). Classification and decomposition of quantum Markov semigroup. *Probab. Theory Related Fields* **134** 603–623. MR2214906
[61] WISEMAN, H. M. (1996). Quantum trajectories and quantum measurement theory. *Quantum Semiclass. Optics* **8** 205–222. MR1374517
[62] ZEIDLER, E. (1990). *Nonlinear Functional Analysis and Its Applications.* II/A. *Linear Monotone Operators*. Springer, New York. MR1033497

LABORATORIO DE ANÁLISIS ESTOCÁSTICO
AND
CENTRO DE ÓPTICA E INFORMACIÓN CUÁNTICA
DEPARTAMENTO DE INGENIERÍA MATEMÁTICA
FACULTAD DE CIENCIAS FÍSICAS Y MATEMÁTICAS
UNIVERSIDAD DE CONCEPCIÓN
CASILLA 160 C, CONCEPCIÓN
CHILE
E-MAIL: cmora@ing-mat.udec.cl

LABORATORIO DE ANÁLISIS ESTOCÁSTICO
DEPARTAMENTO DE MATEMÁTICA
FACULTAD DE MATEMÁTICAS
PONTIFICIA UNIVERSIDAD CATÓLICA DE CHILE
CASILLA 306, SANTIAGO 22
CHILE
E-MAIL: rrebolle@puc.cl